\documentclass[12pt]{article}

\usepackage[latin1]{inputenc}
\usepackage{amsmath,amssymb}
\usepackage{amsmath,accents}
\usepackage{alltt,graphics,graphpap,color}
\usepackage[dvips]{graphicx}
\usepackage{amsfonts}
\usepackage{amscd}
\usepackage{mathrsfs}
\usepackage{hyperref}
\usepackage{float}
\usepackage{subcaption} 

\DeclareMathAlphabet{\mathpzc}{OT1}{pzc}{m}{it}

\newtheorem{Theorem}{Theorem}[section]

\newtheorem{Proposition}[Theorem]{Proposition}
\newtheorem{Lemma}[Theorem]{Lemma}

\newtheorem{Remark}[Theorem]{Remark}

\newcommand{\mm}{{\cal M}}

\newcommand{\proof}{\emph{Proof. }}
\newcommand{\cvd}{\hfill$\square$ \bigskip}
\newcommand{\hhh}{Nil}
\newcommand{\sql}{\sqrt{\lambda}}
\newcommand{\ccc}{\cos(v_2)}

\newcommand{\sss}{\sin(v_2)}

\newcommand{\nuv}{\vec{\nu}}
\newcommand{\tf}{\tilde f}
\newcommand{\tfz}{\tilde f(z)}

\begin{document}

\def\qed{\hbox{\hskip 6pt\vrule width6pt height7pt
depth1pt  \hskip1pt}\bigskip}


\title{Invariant translators of the Heisenberg group}

\author{\sc Giuseppe Pipoli}
\date{}

\maketitle

{\small \noindent {\bf Abstract:} We classify all the translating solitons to the mean curvature flow in the three-dimensional  Heisenberg group that are invariant under the action of some one-parameter group of isometries of the ambient manifold. The problem is solved considering any canonical deformation of the standard Riemannian metric of the Heisenberg group. We highlight similarities and differences with the analogous Euclidean translators: we mention in particular that we describe the analogous of the tilted grim reaper cylinders, of the bowl solution and of translating catenoids, but some of them are not convex in contrast with a recent result of Spruck and Xiao \cite{SX} in the Euclidean space. Moreover we also prove some negative results. Finally we study the convergence of these surfaces as the ambient metric converges to the standard sub-Riemannian metric on the Heisenberg group.} 

\medskip


\bigskip

\section{Introduction}

A hypersurface $\mm$ in a given ambient manifold $(\overline{\mm},\bar g)$ is said a \emph{soliton to the mean curvature flow} if there is a one-parameter group $G=\left\{\varphi_t\left|\ t\in \mathbb R\right.\right\}$ of isometries of $\bar g$ such that the evolution by mean curvature flow starting from $\mm$ is given at any time $t$ by
$$
\mm_t=\varphi_t\left(\mm\right).
$$
Let $V$ be the Killing vector field associated to the group $G$, it turns out that the property of being a soliton can be translated in a prescribed mean curvature problem:
\begin{equation}\label{trasl}
H=\bar g(\nuv,V),
\end{equation}
where $H$ is the mean curvature of $\mm$ and $\nuv$ its unit normal vector field. A proof can be read for example in \cite{HuSm}.

Special solitons are the \emph{translators} in the Euclidean space: $V$ is a constant vector field and therefore $\mm$ evolves translating with constant speed in the direction of $V$.  They have a deep meaning in the analysis of the singularities of the flow: Huisken and Sinestrari \cite{HuSi} proved that convex translators appear as blow-up of type II singularities of mean convex curvature flow. Some interesting examples in this class are known. In the plane there is only one (up to isometries): it is the well known \emph{grim reaper} of equation $y=-\log(\cos(x))$. For higher dimension there is much more freedom: the cylinder generated by the grim reaper is the simplest example, but it can be also deformed in a proper way giving a one-parameter family of tilted grim reaper cylinders: see \cite{BLT} or \cite{HIMW} for an exhaustive description. Clutterbuck, Schn\"urer and Schulze \cite{CSS} described all the rotationally invariant translators: there is only one that is a complete graph, it is called the \emph{bowl solution}, and a one-parameter family of surfaces with two ends. The latter are often called the \emph{translating catenoids}. Very recently Bourni, Langford and Tinaglia \cite{BLT} found a class of graphical not symmetric translators. The biggest achievement in this subject probably is the recent classification of graphical translators in $\mathbb R^3$ due to Hoffman, Ilmanen, Martin and White \cite{HIMW}: they proved that tilted grim reapers, the bowl solutions and the examples of \cite{BLT} are the only complete graphs in the Euclidean $3$-space that are translators. A general classification is far from being understood. There are many non graphical examples (those with helicoidal symmetries \cite{Ha} including some not embedded ones and some others with non trivial topology \cite{DPN,Lo,Smi}) and some rigidity results \cite{CM,IR,Mo,MSS}. The literature on this topic is quite large and these lists are far from being complete.  

It is worth mentioning that in literature there are at least two very interesting ways to transform a translator in a minimal surface. The first one is due to Ilmanen \cite{Il}: $\mm$ is a translator for the Euclidean metric $\bar g$ if and only if $\mm$ is minimal with respect to the conformal metric $\tilde g_p=e^{-\frac{2}{n}\bar g(p,V)}\bar g_p$. The second one has been discovered by Smoczyk \cite{Sm} in case of hypersurfaces and generalized by Arezzo and Sun \cite{AS} for higher codimension submanifolds: it gives a correspondence between translators in the Euclidean space $\mathbb R^n$ and minimal submanifolds in $\mathbb R^n\times\mathbb R$ equipped with a suitable warped product metric. The two approaches are very useful because they allow to use tools from the minimal surface theory in the study of translators. Moreover they are flexible enough to be applied to a wide class of curved ambient manifold see for example \cite{Bu,LM}, and to define more general notions of solitons \cite{ALR}.

In the present paper we consider translating solitons in the $3$-dimensional Heisenberg group.  This space is often denoted with $\hhh_3$ and it can be identified with $\mathbb R^3$ equipped with the following group operation $\star$: for any $(x_1,y_1,z_1),(x_2,y_2,z_2)\in\mathbb R^3$ we have
\begin{equation}\label{operation}
(x_1,y_1,z_1)\star(x_2,y_2,z_2)=\left(x_1+x_2,y_1+y_2,z_1+z_2+\frac{x_1y_2-x_2y_1}{2}\right).
\end{equation}
If we denote with $\partial_x,\partial_y,\partial_z$ the usual coordinates vector fields on $\mathbb R^3$, we can define the following left invariant vector fields on $\hhh_3$
$$
X=\partial_x-\frac y2\partial_z,\quad Y=\partial_y+\frac x2\partial_z,\quad Z=\partial_z.
$$
For any $\lambda>0$, let $\bar g_\lambda$ be the Riemannian metric on $\hhh_3$ such that $\left(X,Y,\lambda^{-\frac 12}Z\right)$ is an orthonormal basis. The isometry group of this metric is independent on $\lambda$ and it is generated by \emph{horizontal rotation} (i.e. usual rotations of the first two coordinates) and left multiplication by any point $p$: they are called \emph{left translations} and they are denoted with $L_p$.

A surface $\mm$ in $\hhh_3$ is called a \emph{translator} if its evolution by mean curvature flow is a left translation of $\mm$ along a fixed direction $p\in\hhh_3$.

In the present paper we will classify all the translators invariant under the action of some one-parameter group of isometries of the ambient space. The one-parameter groups of isometries of $\hhh_3$ are described in Theorem 2 of \cite{FMP}. Since $\hhh_3$ is not an isotropic space, contrary to the Euclidean space, the choice of the direction $p$ is fundamental. The most interesting examples appear for surfaces translating in the vertical direction with unit speed, that is $p=(0,0,\lambda^{-\frac 12})$. We call these surfaces \emph{vertical translators}. In this case $V=\lambda^{-\frac 12}{Z}$ and \eqref{trasl} becomes
\begin{equation}\label{vert trasl}
H=\bar{g}_{\lambda}\left(\nuv,\lambda^{-\frac 12}{Z}\right).
\end{equation}

From our point of view $\hhh_3$ is the simplest ambient manifold where the two relations between translators and minimal surfaces mentioned above fail: in fact both of them require that the Killing vector field is the gradient of a function, but this cannot occur in the Heisenberg group essentially because the distribution generated by $X$ and $Y$ is not integrable. 

The symmetric vertical translators of $\hhh_3$ are described in our first main Theorem.

\begin{Theorem}\label{main 1}
Let $\lambda$ be a fixed positive number and let $G$ be a one-parameter group of isometries of $\hhh_3$.
\begin{itemize}
\item [1)] If $G=\left\{L_{(0,0,u)}\left|\ u\in\mathbb R\right.\right\}$ is the group of vertical translation, then $\mm$ is a $G$-invariant vertical translator if and only if it is a vertical plane, in particular it is minimal.
\item [2)] If $G=\left\{L_{(u,0,cu)}\left|\ u\in\mathbb R\right.\right\}$ for some $c\geq 0$, then there exists an unique $G$-invariant vertical translator: it is a complete horizontal graph, it is not convex, its intrinsic sectional curvature has both signs and it is defined on a slab of width 
$$
\Delta_{\lambda,c}=\frac{2}{\sqrt{\lambda}}\sqrt{1+\lambda c^2}\sinh\left(\frac{\pi\sql}{2}\right).
$$
\item [3)] If $G=\left\{\rho_u\left|\ u\in\mathbb R\right.\right\}$ is the group of horizontal rotations, then there exists a one-parameter family of $G$-invariant vertical translators. All of them are embedded. The parameter is the distance from the vertical axis. If such distance is strictly positive, the surface has two ends and each of them is a horizontal graph outside a compact neighborhood of the origin. In analogy with the Euclidean case we call them \emph{translating catenoids}. If the distance is zero, we have the only $G$-invariant vertical translator which is an entire graph. In analogy with the Euclidean case we can call it the \emph{bowl solution}. In every case the ends of these surfaces have the following asymptotic expansion as $r=\sqrt{x^2+y^2}$ goes to infinity:
\begin{equation}\label{0001}
\frac{r^2}{2\sqrt{\lambda}}+\left\{\begin{array}{ccc}
-\frac{4}{\sqrt{\lambda}(4-\lambda)}\log(r) & \text{if} & \lambda<4;\\
-\log^2(r) & \text{if} & \lambda=4;\\
C_0 r^{1-\frac{4}{\lambda}} & \text{if} & \lambda>4,
\end{array}\right.
\end{equation}
for some constant $C_0$ depending on the initial datum.
\item[4)] If $G=\left\{L_{(0,0,cu)}\circ\rho_u\left| u\in\mathbb R\right.\right\}$ is the group of helicoidal motions around the $z$-axis with pitch $c>0$, then there is a one-parameter family of $G$-invariant vertical translators: the generating curve has exactly one point closest to the origin and consist of two properly embedded arms coming from this point which strictly go away to infinity and spiral infinitely many circles around the origin. The two arms can intersect each other infinitely many times producing a not-embedded surface.
\item[5)] Finally let $p=(x_0,y_0,z_0)$ with $(x_0,y_0)\neq(0,0)$ be a generic non-vertical point and let $G=\left\{L_{up}\circ\rho_u\left| u\in\mathbb R\right.\right\}$ be the associated group of helicoidal motions with general axis, then there are no $G$-invariant vertical translators. 
\end{itemize}
In every cases the uniqueness is up to isometries of $\hhh_3$.
\end{Theorem}

The symmetries of the surfaces considered are preserved by the mean curvature flow (see for example \cite{Pi} for a proof) and help reducing the complexity of the problem. From a technical point of view in this way equation \eqref{vert trasl} becomes an ODE on a suitable function. This ODE describes a curve that produces the vertical translator under the action of the group $G$. This strategy has been successfully applied by Figueroa, Mercuri and Pedrosa in \cite{FMP} for the classification of symmetric constant mean curvature surfaces (including minimal) of the Heisenberg group with $\lambda=1$. Their results hold for every $\lambda$ with minor modifications. In the same spirit, more recently the translators with symmetries have been studied in different ambient manifolds: see for example \cite{Bu, KO, LM}.

The groups considered in part $1)$ and $2)$ of Theorem \ref{main 1} exhaust all left translations: in fact applying a horizontal rotation to the whole ambient space it is always possible to bring us back to this situation. The surfaces described in part $2)$ of Theorem \ref{main 1} are the analogous of the tilted grim reaper cylinders of the Euclidean space. Similar to the Euclidean case they are vertical graphs defined on stripes, however they are not symmetric with respect to the plane $x=0$ and they are not convex. A recent result of Spruck and Xiao \cite{SX} says that all graphical translators in the Euclidean space are convex. Our examples show that in a different ambient space with a richer geometry the properties of translators can be less rigid. 

The examples described in part $3)$ of Theorem \ref{main 1} are the analogous of the Euclidean translators of Clutterbuck, Schn\"urer, and Schulze in \cite{CSS}. Similarly to their work we found that each arm has a quadratic growth. 

From a qualitative point of view, the helicoidal vertical translators in $\hhh_3$ of part $4)$ of Theorem \ref{main 1} are very close to the analogous surfaces in the Euclidean space found by Halldorsson: see Theorem 4.1 of \cite{Ha}. However as we will see in section $5$ the richer geometry on our ambient manifold produce a more difficult ODE. An other consequence of the different geometry of the ambient space is the non existence result of part $5)$ of Theorem \ref{main 1}. In $\mathbb R^3$ a helicoid with a general axis is always isometric to a helicoid with vertical axis. In fact it is sufficient to rotate one axis into the other and apply the result of \cite{Ha} to produce translators of the desired type, but this kind of transformation are not isometries of the Heisenberg group.

When $\lambda$ diverges, the metric $\bar g_{\lambda}$ converges to a sub-Riemannian metric. It is interesting to study how the surfaces described in Theorem \ref{main 1} change with $\lambda$ and in particular what happens when $\lambda$ diverges.

\begin{Theorem}\label{main 2}
As $\lambda$ diverges we have:
\begin{itemize}
\item[1)] for any $c\geq 0$ the surface described in part $2)$ of Theorem \ref{main 1} converges, on any strips $\mathbb R\times K$ where $K$ is a compact subset of $\mathbb R$ in the $C^m$-norm for any $m\in\mathbb N$, to the entire graph of equation
$$
z=\frac{xy}{2}+cx,
$$
this surface is minimal in $\left(\hhh_3,\bar g_{\lambda}\right)$ for any $\lambda$ and hence horizontal-minimal in the sub-Riemannian Heisenberg group;
\item[2)] the bowl solution converges uniformly in $\mathbb R^2$ in the $C^m$-norm for any $m\in\mathbb N$ to a horizontal plane that is a minimal surface in $\left(\hhh_3,\bar g_{\lambda}\right)$ for any $\lambda$ and hence horizontal-minimal in the sub-Riemannian Heisenberg group;
\item[3)] the translating catenoids converge, on any cylinder $\mathbb S^1\times K$ where $K$ is a compact subset of $\mathbb R$ in the $C^m$-norm for any $m\in\mathbb N$, to the surface of equation 
$$\sqrt{x^2+y^2}=\frac{\sqrt{4z^2+c^4}}{c},$$
where in this case $c>0$ represents the minimal distance to the vertical axis. This surface is horizontal-minimal in the sub-Riemannian Heisenberg group, but not minimal in any $\left(\hhh_3,\bar g_{\lambda}\right)$.
\end{itemize}
\end{Theorem}

Note that part $3)$ of Theorem \ref{main 2} is coherent with \eqref{0001}: when $\lambda$ diverges the leading term becomes linear.

In order to complete the classification of the translators invariant under some one-parameter group of isometries, the case of a generic direction remains to be dealt with. We recover from a different point of view some of the surfaces described in Theorem \ref{main 1} and we get some new non-existence results. Let $V=a_1X+a_2Y+(a_3+a_1y-a_2x)Z$ with $(a_1,a_2,a_3)\in\mathbb R^3$ and $(a_1,a_2)\neq(0,0)$ be the Killing vector field associated to the group generated by a generic left translation $L_{(a_1,a_2,a_3)}$. Note that in this case it is not possible to require a translation with unit speed because the norm of $V$ is not constant. 

\begin{Theorem}\label{main 3}
Fix a $\lambda>0$ and let $\mm$ be an invariant surface of $\hhh_3$ which evolves translating in the direction of $V$, then:
\begin{itemize} 
\item[1)] if $\mm$ is invariant for vertical translation then $\mm =\gamma\times\mathbb R$, where $\gamma$ is the grim reaper of the Euclidean plane which evolves translating in the direction $(a_1,a_2)$ with constant speed $\sqrt{a_1^2+a_2^2}$, in particular the solution in this case is independ on $\lambda$;
\item[2)] if $\mm$ is invariant for the group generated by a generic translation $L_{(x_0,y_0,z_0)}$, then a solution exists if and only if there is a $c\in\mathbb R$ such that $(a_1,a_2)=c(x_0,y_0)$;  in this case $\mm$ can be viewed as a vertical translator with constant speed $(a_3-cz_0)\sql$ and can be treated as in part $2)$ of Theorem \ref{main 1} after some rescaling;
\item[3)] $\mm$ cannot be invariant for horizontal rotations;
\item[4)] $\mm$ cannot be a helicoid. 
\end{itemize}
\end{Theorem}

The paper is organized as follows. In Section $2$ we collect some preliminaries on the geometry of the Heisenberg group and of its surfaces. In particular in Lemma \ref{2ff} we present all the tools for computing the mean curvature and the normal vector of a generic surface. These quantities are needed to apply the translator equation \eqref{trasl}. Since many of the surfaces that we want to describe are horizontal graphs, we will be more explicit in this special setting. Moreover in Section $2$ we list also some technical results useful in the analysis of the ODE that we are going to study. Section $3$ concerns vertical translators invariant under the group generated by some left translation: part $1)$ and $2)$ of Theorem \ref{main 1} and part $1)$ of Theorem \ref{main 2} are proved.
Section $4$ is devoted to the rotationally invariant vertical translators: part $3)$ of Theorem \ref{main 1} and part $2)$ and $3)$ of Theorem \ref{main 2} are proved. Helicoidal vertical translators are studied in Section $5$ proving part $4)$ and $5)$ of Theorem \ref{main 1}. Finally in Section $6$ we conclude the classification with the proof of Theorem \ref{main 3}.

\section{Preliminaries}
\subsection{Geometry of $\hhh_3$}

For any $\lambda>0$ let $\bar g_{\lambda}$ be the metric on $\hhh_3$ such that the basis $(X,Y,\lambda^{-\frac 12}Z)$ is orthonormal. We call \emph{horizontal distribution} the distribution $\mathcal H$ generated by $X$ and $Y$. We note that $\left[\mathcal H,\mathcal H\right]=\mathbb RZ$, hence $\mathcal H$ is not integrable and a metric tensor defined only on it and not on the whole tangent space is enough to define a distance. This kind of metrics are called \emph{sub-Riemannian metrics}. When $\lambda$ diverges, the family of Riemannian metric $\bar{g}_{\lambda}$ converges to the \emph{standard sub-Riemannian metric} on $\hhh_3$: that one such that $X$ and $Y$ are orthonormal. The associated distance is called the \emph{Carnot-Carath\'eodory distance}. A very good monograph dedicated to this subject is \cite{CDPT}.

For any finite $\lambda$, the geometry of $\bar g_{\lambda}$ is well known. The Levi-Civita connection with respect to $(X,Y,\lambda^{-\frac 12}Z)$ is:

\begin{equation}\label{Levi Civita}
\bar\nabla = \frac{1}{2}\left(\begin{array}{ccc}
0 &  Z & -\lambda^{\frac 12}Y\\
-Z & 0 & \lambda^{\frac 12}X\\
-\lambda^{\frac 12}Y & \lambda^{\frac 12}X& 0
\end{array}\right)
\end{equation}
The sectional curvature of the metric $\bar g_{\lambda}$ is 
$$
\bar K(X\wedge Y)=-\frac 34\lambda,\quad\bar K(X\wedge Z)=\bar K(Y\wedge Z)=\frac 14\lambda.
$$
For a general tangent plane we can always find a special orthonormal basis of that plane of the type $\left(V_1,V_2=aU+b\lambda^{-\frac 12}Z\right)$, with $V_1, U\in\mathcal H$, $\bar g_{\lambda}(V_1,U)=0$ and $a,b\in\mathbb R$. In this case 
\begin{equation}\label{sec curv}
\bar K(V_1\wedge V_2)=\frac{\lambda}{4}\left(b^2-3a^2\right).
\end{equation}
The group of isometry $Iso(\hhh_3)$ is independent on $\lambda$: it is generated by left translation denoted by $L_p$, that is the multiplication on the left by a fixed point $p\in\hhh_3$, and horizontal rotations, that is usual rotation of the first two components. In the following we denote with $\rho_u$ those who rotate the space with an angle $u$ preserving the orientation:
$$
\rho_u:(x,y,z)\mapsto\left(x\cos(u)-y\sin(u),x\sin(u)+y\cos(u),z\right).
$$
In particular $Iso(\hhh_3)$ has dimension $4$ that is the maximum that we can expect from a $3$-dimensional manifold that is not a space form. A basis of the Lie algebra of the Killing vector field is:
$$\begin{array}{rclcrcl}
F_1  &=&  \partial_x+\frac y2\partial_z = X+yZ; && F_2  &=&  \partial_y-\frac x2\partial_z = Y-xZ;\\
F_3  &=&  \partial_z = Z; && F_4  &=&  -y\partial_x+x\partial_y = -yX+xY-\frac{(x^2+y^2)}{2}Z.
\end{array}
$$
The first three are generated by left translation and the latter by the rotations. Note that $F_3$ is the only one with constant norm. The closed one-dimensional subgroup of $Iso(\hhh_3)$ are generated by linear combination with constant coefficients of the four Killing vector field above, as shown in Theorem 2 of \cite{FMP}.

\subsection{Surfaces of $\hhh_3$}
The following result will be used quite often in this paper because it provides us with the tools we need to compute the quantities that appear in the translating equation \eqref{trasl}.

\begin{Lemma}\label{2ff}
Let $\mm$ be a surface of $\left(\hhh_3,\bar g_{\lambda}\right)$ parametrized by 
$$
\mm(v_1,v_2)=\left(x(v_1,v_2),y(v_1,v_2),z(v_1,v_2)\right).
$$
Then:
\begin{itemize}
\item[1)] a basis of the tangent space of $\mm$ is given by the vector fields $V_1, V_2$ where for every $i=1,2$  $V_i=a_iX+b_iY+c_iZ$ and
$$
a_i=\frac{\partial x}{\partial v_i},\ b_i=\frac{\partial y}{\partial v_i},\ c_i=\frac{\partial z}{\partial v_i}+\frac{1}{2}\left(y\frac{\partial x}{\partial v_i}-x\frac{\partial y}{\partial v_i}\right);
$$
\item[2)] for every $i,j=1,2$ we have
\begin{eqnarray*}
\bar{\nabla}_{V_i}V_j & = & \left(\frac{\partial a_j}{\partial v_i}+\frac{{\lambda}}{2}\left(b_jc_i+b_ic_j\right)\right)X\\
&&+\left(\frac{\partial b_j}{\partial v_i}-\frac{{\lambda}}{2}\left(a_jc_i+a_ic_j\right)\right)Y\\
&&+\left(\frac{\partial c_j}{\partial v_i}+\frac{1}{2}\left(b_ja_i-b_ia_j\right)\right)Z.
\end{eqnarray*}
\end{itemize}
Moreover, if $\mm$ is an horizontal graph, i.e. $z=u(x,y)$, using the coordinates $v_1=x$ and $v_2=y$, let us define
$$
\alpha=u_x+\frac y2,\quad \beta=u_y-\frac x2,
$$
then:
\begin{itemize}
\item[3)] $V_1=X+\alpha Z$, $V_2=Y+\beta Z$ and a unit normal vector field is $$\nuv=\frac{1}{\sqrt{\lambda+\lambda^2\alpha^2+\lambda^2\beta^2}}\left(-\alpha\lambda X-\beta\lambda Y+Z\right);$$
\item[4)] with respect $V_1,V_2$ induced metric $g$ and the second fundamental form $A$ and the mean curvature $H$ of $\mm$ are
$$
g=\left(\begin{array}{cc}
1+\alpha^2\lambda & \alpha\beta\lambda\\
\alpha\beta\lambda & 1+\beta^2\lambda
\end{array}\right),
$$
$$
A=\frac{\lambda}{\sqrt{\lambda+\lambda^2\alpha^2+\lambda^2\beta^2}}\left(\begin{array}{cc}
u_{xx}+\alpha\beta\lambda & u_{xy}+\frac{\lambda}{2}\left(\beta^2-\alpha^2\right)\\
u_{xy}+\frac{\lambda}{2}\left(\beta^2-\alpha^2\right) & u_{yy}-\alpha\beta\lambda
\end{array}\right)
$$
and
$$
H=\frac{\sqrt{\lambda}}{\left(1+\alpha^2\lambda+\beta^2\lambda\right)^{\frac 32}}\left(u_{xx}(1+\beta^2\lambda)+u_{yy}(1+\alpha^2\lambda)-2u_{xy}\alpha\beta\lambda\right).
$$
\end{itemize}
\end{Lemma}
\proof
\begin{itemize}
\item[1)] The vector fields 
$$V_i=\frac{\partial\mm}{\partial v_i}=\frac{\partial x}{\partial v_i}\frac{\partial}{\partial x}+\frac{\partial y}{\partial v_i}\frac{\partial}{\partial y}+\frac{\partial z}{\partial v_i}\frac{\partial}{\partial z}
$$
clearly span the tangent space of $\mm$. Writing them in the basis $(X,Y,Z)$ we have the thesis.
\item[2)] Fix $i$ and $j$, then, using the linearity and the Leibniz rule of the Levi-Civita connection we have:
\begin{eqnarray*}
\bar{\nabla}_{V_i}V_j & = & \frac{\partial a_j}{\partial v_i}X+\frac{\partial b_j}{\partial v_i}Y+\frac{\partial c_j}{\partial v_i}Z\\
 & & +a_j\left(a_i\bar{\nabla}_XX+b_i\bar{\nabla}_YX+c_i\bar{\nabla}_ZX\right)\\
& & +b_j\left(a_i\bar{\nabla}_XY+b_i\bar{\nabla}_YY+c_i\bar{\nabla}_ZY\right)\\
& & +c_j\left(a_i\bar{\nabla}_XZ+b_i\bar{\nabla}_YZ+c_i\bar{\nabla}_ZZ\right).
\end{eqnarray*}
Applying the explicit expression of $\bar{\nabla}$ \eqref{Levi Civita} we get the thesis.
\item[3)] In this special setting we have that $a_i=\delta_{i1}$, $b_i=\delta_{i2}$, $c_1=\alpha$ and $c_2=\beta$. Moreover it easy to check that
$$
\bar g_{\lambda}(\nuv,\nuv)=1,\quad\bar g_{\lambda}(\nuv,V_i)=0.
$$
\item[4)] Recalling that $g_{ij}=\bar{g}_{\lambda}(V_i,V_j)$, $h_{ij}=A(V_i,V_j)=\bar{g}_{\lambda}\left(\bar{\nabla}_{V_i}V_j,\nuv\right)$ and that $H=h_{ij}g^{ij}$, the thesis follows by parts $2)$  and $3)$ of the present Lemma.\cvd
\end{itemize}

We conclude this section recalling some basic definitions about the geometry of submanifolds in the sub-Riemannian setting. We refer to \cite{CDPT} for all the details. Given a $C^2$-surface $\mm$ in $\hhh_3$, we call \emph{characteristic points} of $\mm$ the points $p\in\mm$ such that the tangent space in $p$ coincides with $\mathcal H_p$. The other points of the surface are called \emph{noncharacteristic points}. For any given $\lambda>0$, let $H(\lambda)$ be the mean curvature of $\mm$ with respect to $\bar g_{\lambda}$. One can prove that in any noncharacteristic point the $\lim_{\lambda\rightarrow\infty}H(\lambda)$ exists and is finite. This limit is called the \emph{horizontal mean curvature} of $\mm$. In analogy with the Riemannian case, a surface is said \emph{horizontal-minimal} if its horizontal mean curvature vanishes along its noncharacteristic points.

\subsection{Some technical results}

We introduce some technical analytic results used several times in the present paper. The proof of the first one is a straightforward calculation, so we omit it.

\begin{Lemma}\label{ode speciale}
For any $a,b,c\in\mathbb R$ with $b>0$, let $y$ be a solution of the following ODE
$$
y'=-\frac{b}{x}(y-a)+\frac{c}{x^2},
$$ 
then 
$$
y=\left\{\begin{array}{ccc}
\displaystyle{a+\frac dx+c\frac{\log(x)}{x} }& \text{if} & b=1;\\
\\
\displaystyle{a+\frac{d}{x^b}+\frac{c}{b-1}\frac{1}{x}}& \text{if} & b\neq 1;
\end{array}\right.
$$
for some constant $d$. In both case $y$ is defined for any $x$ bigger than some $x_0>0$ and
$$
\lim_{x\rightarrow\infty}y(x)=a.
$$
\end{Lemma}

The second one will be useful when we will investigate the convergence for $\lambda$ going to infinity and proving Theorem \ref{main 2}. Here and in the following we denote with $f^{(0)}=f$ and for any $m\geq 1$, with $f^{(m)}$ the $m$-th derivative of $f$.

\begin{Lemma}\label{derivate}
Let $f,g$ two smooth real functions, then for every $m\in\mathbb N$ we have:
\begin{itemize}
\item[1)] $\left(fg\right)^{(m)}=\sum_{j=0}^m {m\choose j} f^{(j)}g^{(m-j)}$;
\item[2)] there exists a polynomial $P_m$ of degree $m$ with integer coefficients and $m+1$ variables such that 
$$
\left(\frac 1 f\right)^{(m)}=\frac{P_s(f^{(0)},f^{(1)},\dots,f^{(m)})}{f^{m+1}}.
$$
\end{itemize}
\end{Lemma}
\proof
\begin{itemize}
\item[1)] This equality can be easily proved by induction on $m$.
\item[2)] We proceed by induction on $m$. If $m=0$  or $m=1$ the statement is trivial. Fix an $m\geq 2$ and suppose that the result is true for $m$. We have:
\begin{eqnarray*}
\left(\frac 1f\right)^{(m+1)} & = & \left(\frac {P_m(f,\cdots,f^{(m)})}{f^{m+1}}\right)^{(1)}\\
& = & \frac{f\sum_{j=0}^m f^{(j+1)}\frac{\partial P_m}{\partial f^{(j)}}-(m+1)f^{(1)}P_m}{f^{m+2}}.
\end{eqnarray*}
Defining $P_{m+1}(f^{(0)},\dots,f^{(m)},f^{(m+1)})=f\sum_{j=0}^m f^{(j+1)}\frac{\partial P_m}{\partial f^{(j)}}-(m+1)f^{(1)}P_m$ we get the thesis. \cvd
\end{itemize}

\section{Surfaces invariant under translations}
In this section we will describe the vertical translators invariant under the action of a group generated by a left translation $L_p$ for some $p=(x_0,y_0,z_0)\in\hhh_3$  They are the analogous of the Euclidean tilted grim reaper cylinders. Moreover we will prove also the first part of Theorem \ref{main 2} studying the behavior of such surfaces when $\lambda$ diverges. 

Note that, applying a suitable horizontal rotation to the whole ambient space, we can always assume that $y_0=0$. This fact simplify the computations and helps in avoiding repetitions up to isometries of $\hhh_3$.

The following result proves part $1)$ and $2)$ of Theorem \ref{main 1}.
\begin{Theorem}\label{teor trasl 2} 
Fix a $\lambda >0$, then we have:
\begin{itemize}
\item[1)] a vertical translator $\mm$ is invariant under vertical translations if and only if it is a vertical plane, in particular it is minimal;
\item[2)] for every $c\geq 0$ let $G=\left\{L_{(u,0,cu)}\left|\ u\in\mathbb R\right. \right\}$ be the group of translations with slope $c$, then there exists a $G$-invariant vertical translator $\mm^{\lambda,c}$: it is a complete horizontal graphs, it is not convex, its intrinsic sectional curvature has both signs and it is defined on a slab of width 
$$
\Delta_{\lambda,c}=\frac{2}{\sqrt{\lambda}}\left(1+\lambda c^2\right)^{\frac 12}\sinh\left(\frac{\pi\sql}{2}\right).
$$
Moreover $\mm^{\lambda,c_1}$ and $\mm^{\lambda,c_2}$ are not equivalent under an isometry of the ambient manifold if $c_1\neq c_2$. 
\end{itemize}
\end{Theorem}
\proof
\begin{itemize}
\item[1)]  If $\mm$ is invariant under vertical translations, then $Z$ is a tangent vector field at any point of $\mm$. Therefore the translating property \eqref{vert trasl} says that $\mm$ is minimal. By \cite{FMP} we know that, in this particular case, $\mm$ is a vertical plane.
\item[2)]
Fix a $\lambda>0$ and a $c\geq 0$, since we want to find a $G$-invariant surface $\mm^{\lambda,c}$, we look for a planar curve $(\mu(y),\gamma(y))$ such that $\mm^{\lambda,c}$ is parametrized by
$$
(x,0,cx)\star(0,\mu(y),\gamma(y))=\left(x,\mu(y),\frac{x\mu(y)}{2}+cx+\gamma(y)\right)
$$
is a vertical translator. We can start simplifying the problem restricting our attention to the case of a graph, i.e. $\mu(y)=y$. We will see in a while that it is not actually a restriction because we will prove that such graph is a complete curve. With respect to the notation introduced in Lemma \ref{2ff}, in our case we have
$$
\alpha=y+c,\quad \beta=\gamma'(y).
$$ By Lemma \ref{2ff} part $4)$ we have that
\begin{equation}\label{mc graph lambda}
H=\frac{\sqrt{\lambda}}{\left(1+\alpha^2\lambda+\beta^2\lambda\right)^{\frac 32}}\left(\gamma''(y)(1+\alpha^2\lambda)-\alpha\beta\lambda\right).
\end{equation}
Imposing the equation of translators \eqref{vert trasl} we have that $\gamma$ is the solution of the following Cauchy problem:
\begin{equation}\label{Cauchy lambda}
\left\{\begin{array}{l}
\gamma''=\frac{1}{\sql}+\frac{\sql (\gamma')^2+\lambda(y+c)\gamma'}{1+\lambda (y+c)^2},\\
\gamma(0)=\gamma'(0)=0.
\end{array}\right.
\end{equation}
Luckily  we are able to find the explicit expression for $\gamma'$ in term of elementary functions:
\begin{equation}\label{gamma'}
\gamma'(y)=\sqrt{(y+c)^2+\frac{1}{\lambda}}\tan\left(\frac{1}{\sql}\sinh^{-1}(\sql (y+c))-\frac{1}{\sql}\sinh^{-1}(\sql c)\right).
\end{equation}

Such function is defined on the interval $(a_{\lambda,c},b_{\lambda,c})$ where
\begin{eqnarray*}
a_{\lambda,c} & = & \frac{1}{\sql}\sinh\left(\sinh^{-1}(\sql c)-\frac{\pi}{2}\sql\right)-c;\\
b_{\lambda,c} & = & \frac{1}{\sql}\sinh\left(\sinh^{-1}(\sql c)+\frac{\pi}{2}\sql\right)-c.
\end{eqnarray*}

The width of this interval is 
$$
\Delta_{\lambda,c}=b_{\lambda,c}-a_{\lambda,c}=\frac{2}{\sqrt{\lambda}}\sqrt{1+\lambda c^2}\sinh\left(\frac{\pi}{2}\sql\right).
$$
It is easy to see that $y\gamma'(y)>0$ for every $y\neq 0$, then $\gamma$ has a global minimum in $y=0$. Moreover $\gamma$ blows up when $y$ approaches the extrema of the interval of definition. Therefore $\gamma$ is enough to generate a complete surface and therefore $\mm^{\lambda,c}$ is a horizontal graph. More precisely we have  
$$
\begin{array}{rclcc}
\gamma(y) & \approx & -\frac{\left(1+\lambda(c+b_{\lambda,c})^2\right)}{\sql}\log\left(b_{\lambda,c}-y\right) & \text{as} & y\approx b_{\lambda,c};\\
\gamma(y) & \approx & -\frac{\left(1+\lambda(c+a_{\lambda,c})^2\right)}{\sql}\log\left(y-a_{\lambda,c}\right) & \text{as} & y\approx a_{\lambda,c}.
\end{array}
$$
In fact when $y\rightarrow b_{\lambda,c}$ we get:
$$
\frac{\frac{1}{\sql}\sinh^{-1}(\sql (y+c))-\frac{1}{\sql}\sinh^{-1}(\sql c)-\frac{\pi}{2}}{y-b_{\lambda,c}}\rightarrow \frac{1}{\sqrt{1+\lambda\left(c+b_{\lambda,c}\right)^2}}.
$$
Hence
\begin{eqnarray*}
\gamma'(y) & \approx & \frac{\sqrt{1+\left(c+b_{\lambda,c}\right)^2}}{\sql}\tan\left(\frac{\pi}{2}+\frac{y-b_{\lambda,c}}{\sqrt{1+\left(c+b_{\lambda,c}\right)^2}}\right)\approx  -\frac{1+\lambda\left(c+b_{\lambda,c}\right)^2}{\sql\left(y-b_{\lambda,c}\right)}.
\end{eqnarray*}
Integrating this estimate we have the behavior of $\gamma$ near $b_{\lambda,c}$. In an analogous way we can compute the behavior in $a_{\lambda,c}$.

Specifying what found in Lemma \ref{2ff} part $4)$ in our case we have that
$$
A=(1+\lambda\alpha^2+\lambda\beta^2)^{-\frac 12}\left(\begin{array}{cc}
\lambda\alpha\beta &\frac{\lambda\beta^2-\lambda\alpha^2+1}{2}\\
\frac{\lambda\beta^2-\lambda\alpha^2+1}{2} & \gamma''(y)-\lambda\alpha\beta
\end{array}\right).
$$
The inverse of the induced metric is
$$
g^{-1}=(1+\lambda\alpha^2+\lambda\beta^2)^{-1}\left(\begin{array}{cc}
1+\lambda\beta^2&-\lambda\alpha\beta\\
-\lambda\alpha\beta & 1+\lambda\alpha^2
\end{array}\right).
$$
The Gaussian curvature is, after some computations, 
\begin{equation}\label{gaussian}
\det(Ag^{-1})=-\frac{\lambda\left(\lambda^2\alpha^4+\lambda^2\alpha^2\beta^2-2\lambda\alpha^2-4\sql\alpha\beta+\lambda\beta^2+1\right)}{4(1+\lambda\alpha^2)(1+\lambda\alpha^2+\lambda\beta^2)},
\end{equation}
where, by \eqref{Cauchy lambda}, we used $\gamma''=\frac{1}{\sql}+\frac{\sql\beta^2+\lambda\alpha\beta}{1+\lambda\alpha^2}$. In the special case $y=0$ we have that $\alpha=c$ and $\beta=0$, therefore \eqref{gaussian} becomes:
$$
\det(Ag^{-1})=-\frac{\lambda(1-\lambda c^2)^2}{4(1+\lambda c^2)^2}<0
$$
if $\lambda c^2\neq 1$. Hence in the general case $\mm^{\lambda,c}$ is not convex. Now we want to compute the intrinsic sectional curvature of this surface. From Lemma \ref{2ff} part $1)$ we have that a basis of the tangent space of $\mm^{\lambda,c}$ is $V_1=X+\alpha Z$ and $V_2=Y+\beta Z$. From this basis we can find an orthonormal basis $(\tilde V_1,\tilde V_2)$:
\begin{equation*}
\tilde V_1=\frac{\beta X-\alpha Y}{\sqrt{\alpha^2+\beta^2}},\quad
\tilde V_2=\frac{\tilde U}{\sqrt{1+\lambda\alpha^2+\lambda\beta^2}}+\sqrt{\frac{\lambda\alpha^2+\lambda\beta^2}{1+\lambda\alpha^2+\lambda\beta^2}}\frac{Z}{\sql},
\end{equation*}
where $\tilde V_1,\ \tilde U=\frac{\alpha X+\beta Y}{\sqrt{\alpha^2+\beta^2}}\in\mathcal H$ and they are orthonormal. By \eqref{sec curv} 
$$
\bar K(T\mm^{\lambda,c})=\frac{\lambda\left(\lambda\alpha^2+\lambda\beta^2-3\right)}{4\left(1+\lambda\alpha^2+\lambda\beta^2\right)}.
$$
Applying the Gauss equation and \eqref{gaussian} we have that the intrinsic sectional curvature is
\begin{equation}\label{intrinsic sec}
K(T\mm^{\lambda,c})=\bar K(T\mm^{\lambda,c})+\det(Ag^{-1})=\frac{\lambda\left(\sql\alpha\beta-1\right)}{\left(1+\lambda\alpha^2\right)\left(1+\lambda\alpha^2+\lambda\beta^2\right)}
\end{equation}
If $y=0$, then $\beta=0$ and $K(T\mm^{\lambda,c})<0$. On the other hand, if $b_{\lambda,c}+c>0$, letting $y\rightarrow b_{\lambda,c}$, we have that $\alpha\beta\rightarrow +\infty$ and then $K(T\mm^{\lambda,c})>0$ far enough from $y=0$. If $b_{\lambda,c}+c\leq 0$ we can get the same result when $y\rightarrow a_{\lambda,c}$. Therefore in any case $K(T\mm^{\lambda,c})$ assume both signs. Note that, when $y$ approaches the extrema of the interval $K(T\mm^{\lambda,c})$ tends to zero. 

Finally we show that the examples found so far are not equivalent under some isometry of the ambient manifold. Let $c_1, c_2$ be two constants and let $\gamma_1$ (resp $\gamma_2$ ) be the curve solution of \eqref{Cauchy lambda} with parameter $c_1$ (resp. $c_2$). Clearly $\mm^{\lambda,c_1}$ cannot be obtained from $\mm^{\lambda,c_2}$ with an horizontal rotation. Suppose that there is a point $p=(x_0,y_0,z_0)\in\hhh_3$ such that $\mm^{\lambda,c_1}=L_{p}\mm^{\lambda,c_2}$. Then for every $x\in\mathbb R$, $p\star(x,0,c_2x)\in\mm^{\lambda,c_1}$. In particular by \eqref{operation} and the definition of $\mm^{\lambda,c_1}$ we have that
$$
z_0+c_{2}x-\frac{y_{0}x}{2}=\frac{(x+x_0)y_0}{2}+c_1(x+x_0)+\gamma_1(y_0)
$$
holds for every $x$. It follows that $c_2=c_1+y_0$ and $z_0=\frac{x_0y_0}{2}+c_1x_0+\gamma_1(y_0)$. On the other hand for every $y$ where $\mm^{\lambda,c_2}$ is defined $p\star(0,y,\gamma_2(y))\in\mm^{\lambda,c_1}$. Therefore 
$$
\gamma_2(y)=\gamma_1(y+y_0)-\gamma_1(y_0)
$$
holds for every $y$. Deriving this equality, in $y=0$ we get that $\gamma_1'(y_0)=\gamma_2'(0)=0$. Therefore by \eqref{gamma'} we deduce that $y_0=0$. Hence $c_2=c_1$, $z_0=c_1x_0$ and $L_{p}$ is trivial, in the sense that it belongs to the group of isometries of the vertical translator considered. \cvd
\end{itemize}

\begin{figure}[H]
\centering
\includegraphics[width=0.7\textwidth]{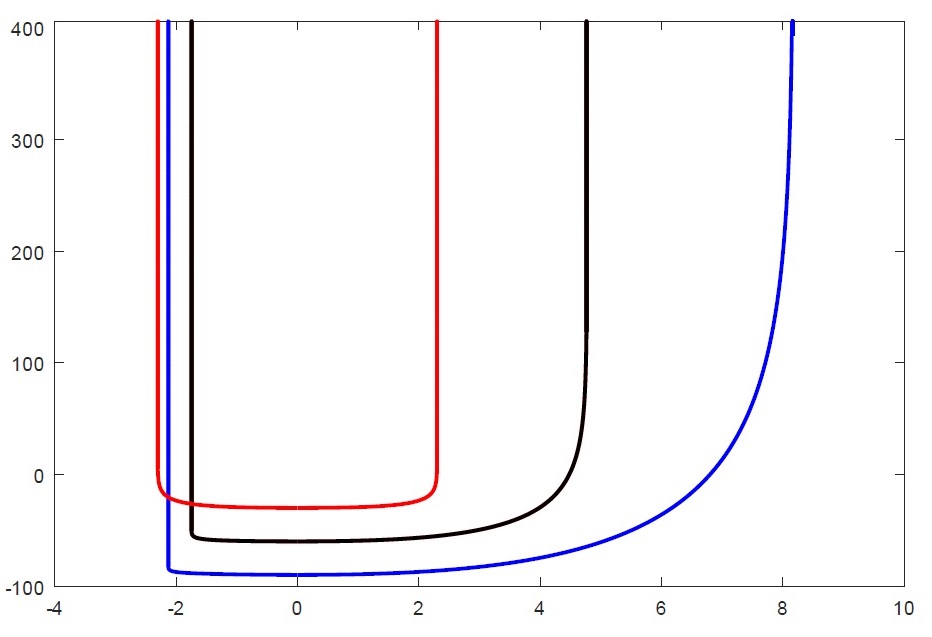}
\caption{Some examples of the generating curve $\gamma$: $\lambda=1$ in any case, while $c=0$ (red curve), $c=1$ (black curve), $c=2$ (blue curve). The graphs have translated in the vertical direction in order to minimize overlapping.}
\end{figure}

\begin{figure}[H]
\centering
\includegraphics[width=0.7\textwidth]{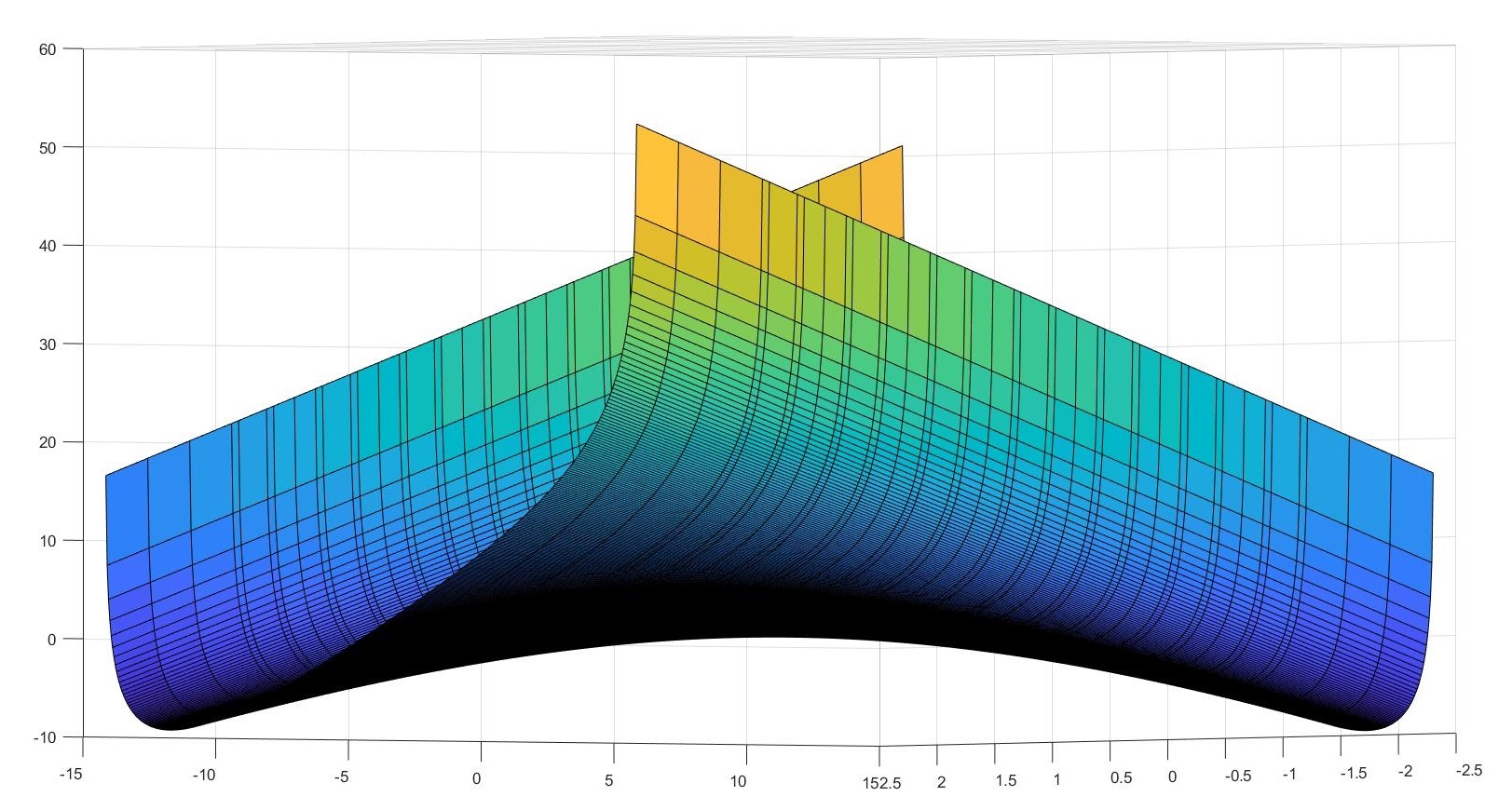}
\caption{A $3d$ rendering of $\mm^{1,0}=\left(x,y,\frac{xy}{2}+\gamma(y)\right)$: one can note the contribution of the quadratic term typical of left translation in $\hhh_3$.}
\end{figure}

\begin{Remark}
The initial condition in \eqref{Cauchy lambda} are not restrictive and allows us to avoid repetitions up to isometries of the ambient space. Since we want to find vertical translators, clearly the height of the surface at time $t=0$ is not important. Hence we can require $\gamma(0)=0$. Moreover one can check that if $\gamma'(0)$ is arbitrary, then we find a solution of \eqref{Cauchy lambda} with a possibly different $c$.
\end{Remark}

Now we want to investigate the behaviour of $\mm^{\lambda,c}$ as $\lambda$ goes to infinity. We start deriving some properties of the limit (if it exists). For any fixed $y$ and $c$, as $\lambda\rightarrow\infty$, the quantity
$$
H=\bar g_{\lambda}\left(\nuv,\frac{Z}{\sql}\right)=\frac{1}{\sqrt{1+\lambda((\gamma'(y))^2+(y+c)^2)}}\rightarrow 
\left\{\begin{array}{ll}
1 & \text{if } y=c=0;\\
0 & \text{otherwise}.
\end{array}\right.
$$
It can be easily checked that $\mm^{\lambda,c}$ has some characteristic points (see Section 2.2 for the definition) if and only if $c=0$ and in this case the characteristic points are those with $y=0$. Therefore, if we have convergence, the limit is a horizontal-minimal surface in the sub-Riemannian Heisenberg group.
Moreover if a limit exists, it is an entire graph, in fact when $\lambda$ diverges we have $
\lim_{\lambda\rightarrow\infty}a_{\lambda,c}=-\infty$ and  $\lim_{\lambda\rightarrow\infty}b_{\lambda,c}=+\infty$.

\begin{Theorem}\label{trasl inf}
Fix any $c\geq 0$, as $\lambda$ diverges, the family of curves $\gamma=\gamma_{\lambda}$ solution of \eqref{Cauchy lambda} converges  to the constant function zero on the compact subsets of $\mathbb R$ in the norm $C^m$ for any $m$. Therefore the family of vertical translators $\mm^{\lambda,c}$ converges in the $C^m$-norm to the surface $\mm^{\infty,c}=\left(x,y,\frac{xy}{2}+cx\right)$. This surface is minimal in $\left(\hhh_3,\bar g_{\lambda}\right)$ for every $\lambda$ and horizontal-minimal in the sub-Riemannian Heisenberg group.
\end{Theorem}


\noindent This Theorem is a direct consequence of the following result.

\begin{Proposition}
For any fixed $c\geq 0$, $m\in\mathbb N$ and compact $K\subset\mathbb R$ there is a positive constant $C=C(c,m,K)$ such that  
$$
\left|\gamma_{\lambda}^{(m)}(y)\right|\leq C\frac{\log(\sql)}{\sql}
$$
holds for every $y\in K$ and any $\lambda$ big enough.
\end{Proposition}
\proof
Fix any $c\geq 0$  and a compact $K$, then if $\lambda$ is big enough we have that $K\subset(a_{\lambda,c},b_{\lambda,c})$. By \eqref{gamma'} we can see that as $\lambda$ diverges
$$
\gamma_{\lambda}'(y)\approx \frac{|y+c|}{\sql}\log\left(\frac{\sqrt{1+\lambda(y+c)^2}+\sql (y+c)}{\sqrt{1+\lambda c^2}+\sql c}\right)
$$
holds. Therefore the thesis follows for $m=1$. Integrating we can recover the case $m=0$ too. The case $m=2$ can be easily derived by \eqref{Cauchy lambda} using the estimate for $m=1$. Fix now an $m\geq 2$ and suppose that the thesis holds for any $j=0,\dots,m$. By Lemma \ref{derivate} end the ODE \eqref{Cauchy lambda} we have:
\begin{eqnarray*}
\gamma_{\lambda}^{(m+1)} & = & \sum_{j=0}^{m-1}{m-1 \choose j} \left(\frac{1}{1+\lambda(y+c)^2}\right)^{(j)}\left(\sql \left(\gamma_{\lambda}'\right)^2+\lambda(y+c)\gamma_{\lambda}'\right)^{(m-1-j)}\\
& = & \sum_{j=0}^{m-1}{m-1 \choose j}\frac{P_j\left(1+\lambda(y+c)^2,\dots,\left(1+\lambda(y+c)^2\right)^{(j)}\right)}{\left(1+\lambda(y+c)^2\right)^{j+1}}\\
&&\phantom{aaaaaa} \cdot\left[\sum_{k=0}^{m-1-j}{m-1-j\choose k}\left(\sql\gamma_{\lambda}^{(k+1)}\gamma_{\lambda}^{(m-j-k)}+\lambda(y+c)^{(k)}\gamma_{\lambda}^{(m-j-k)}\right)\right].
\end{eqnarray*}
Since $\left(\lambda(y+c)^2\right)^{(j)}\approx \lambda$ if $j\leq 2$ and it is zero otherwise, then
$$
\frac{P_j\left(1+\lambda(y+c)^2,\dots,\left(1+\lambda(y+c)^2\right)^{(j)}\right)}{\left(1+\lambda(y+c)^2\right)^{j+1}}\approx\frac{1}{\lambda}.
$$
Finally the term in the square brackets can be estimate with the triangle inequality and the inductive hypothesis getting that it grows at most as $\sql \log(\sql)$. 
\cvd

\begin{Remark} Note that $\gamma=0$ is the the trivial solution of the limit of the Cauchy problem \eqref{Cauchy lambda}:
\begin{equation*}
\left\{\begin{array}{l}
\gamma''=\frac{\gamma'}{y+c},\\
\gamma(0)=\gamma'(0)=0.
\end{array}\right.
\end{equation*}
When $c=0$ the functions $ay^2$ and $ay|y|$ are other solutions of this Cauchy problem for any $a\in\mathbb R$ that do not appears as limit of translators.
\end{Remark}

\begin{figure}[H]
\centering
\includegraphics[width=0.7\textwidth]{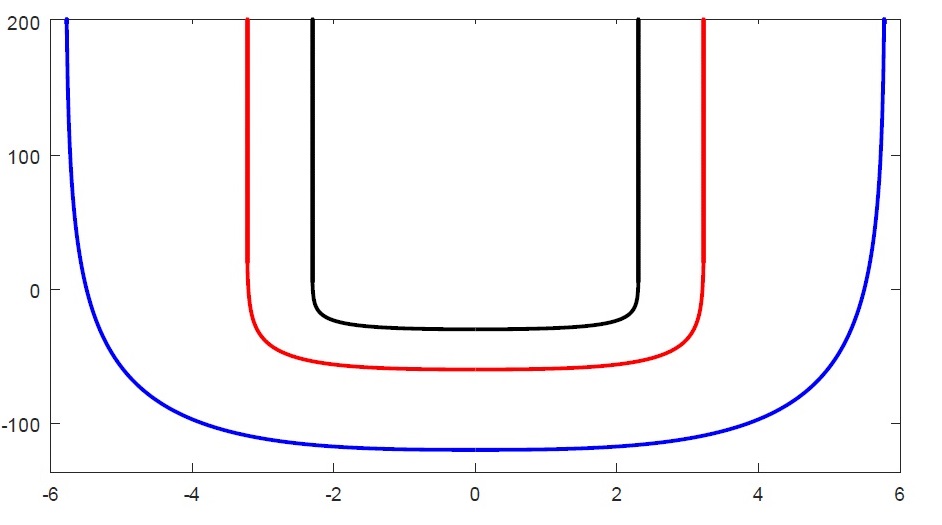}
\caption{Some examples of the generating curve $\gamma$: $c=0$ in any case and $\lambda=1$ (black curve), $\lambda=2$ (red curve) and $\lambda=4$ (blue curve). The curves have been translated vertically in order to avoid overlapping.}
\end{figure}

\section{Surfaces of revolution}
In this section we will prove part $3)$ of Theorem \ref{main 1} producing rotationally invariant vertical translators. The qualitative result is very close to the analogous surfaces described in \cite{CSS}. The main difference is that in our case the asymptotic behavior of the solution depends on $\lambda$ (see Lemma \ref{rotazione stime}). Moreover we will finish the proof of Theorem \ref{main 2} studying the convergence of these translators when $\lambda$ diverges.

\begin{Theorem}\label{teor rot}
For any $\lambda>0$ there is, up to isometry of $\hhh_3$, a one-parameter family of rotationally invariant vertical translators. All of them are embedded. The parameter is the distance from the vertical axis. If such distance is positive, the surface has two ends.  Moreover there is only one vertical translator which is a rotationally invariant entire graph: it can be thought as the surface at distance zero from the vertical axis. 
\end{Theorem}
\proof
We start considering a graph of revolution in $\left(\hhh_3,\bar g_{\lambda}\right)$: let $r=\sqrt{x^2+y^2}$ be the radial coordinate, then there is a function $\varphi:I\subset\mathbb R^+\rightarrow\mathbb R$ such that our surface is the graph $z=u(x,y)=\varphi(r)$. Like in the previous section we start computing the mean curvature of this surface with the help of Lemma \ref{2ff}. Since 
$$
r_x=\frac xr,\quad r_y=\frac yr
$$
we have:
$$
\begin{array}{l}
\displaystyle{u_x = \frac xr\varphi',\quad\quad u_y = \frac yr\varphi',}\\
\displaystyle{u_{xx} = \frac{x^2}{r^2}\varphi''+\frac{1}{r}\varphi'\left(1-\frac{x^2}{r^2}\right),}\\
\displaystyle{u_{yy} = \frac{y^2}{r^2}\varphi''+\frac{1}{r}\varphi'\left(1-\frac{y^2}{r^2}\right),}\\
\displaystyle{u_{xy} = u_{yx} = \frac{xy}{r^2}\varphi''-\frac{xy}{r^3}\varphi',}\\
\displaystyle{\alpha = \frac xr\varphi'+\frac y2,\quad\quad\beta = \frac yr\varphi'-\frac x2,}
\end{array}
$$
where the indices denote the partial derivatives of the functions. By part $4)$ of Lemma \ref{2ff}, after some computations we have:
\begin{equation}\label{H rot}
H=\frac{\sqrt{\lambda}}{\left(1+\lambda\left(\left(\varphi'\right)^2+\frac{r^2}{4}\right)\right)^{\frac 32}}\left[\varphi''\left(1+\frac{\lambda}{4}r^2\right)+\frac{1}{r}\varphi'\left(1+\lambda\left(\varphi'\right)^2\right)\right].
\end{equation}
Imposing the vertical translators equation \eqref{vert trasl}, after some algebraic manipulation we have:
\begin{equation}\label{ode rotazione}
\varphi'' = \frac{1}{\sqrt{\lambda}}+\frac{4\varphi'}{r\left(4+\lambda r^2\right)}\left(\sqrt{\lambda} r \varphi' -1-\lambda \left(\varphi'\right)^2\right).
\end{equation}

As in the Euclidean space \cite{CSS}, we see that there exists only one (up to a vertical translation) of such surfaces that is an entire graphs, i.e. the function $\varphi$ is defined globally for any $r\geq 0$: it is sufficient to require that 
$$
\lim_{r\rightarrow 0}\frac{\varphi'}{r}=\frac{1}{2\sqrt{\lambda}}.
$$
In all the other cases, the function $\varphi$ is defined for $r\geq r_0$ for some positive $r_0$.

Following the idea of \cite{CSS} we can glue two of these graphical solutions for producing a complete rotationally invariant vertical translators. The resulting surface has therefore the topology of the cylinder. In order to do it we need to find a vertical translators of the following type:
\begin{equation}\label{sup rot2}
\mathcal M(v_1,v_2)=\left(f(v_1)\cos(v_2),f(v_1)\sin(v_2),v_1\right)
\end{equation}
for some real smooth function $f$. As before we use Lemma \ref{2ff}. For simplicity of notation, we will denote $f(v_1)$ simply by $f$ and analogously for its derivatives.  We have
$$
\begin{array}{lll}
a_1=f'\ccc, & b_1=f'\sss, & c_1=1,\\
a_2=-f\sss, & b_2=f\ccc, & c_2=-\frac{f^2}{2}.\\
\end{array}
$$
Then a tangent basis is 
$$
V_1=f'\ccc X+f'\sss Y+Z,\quad V_2=-f\sss X+f\ccc Y-\frac{f^2}{2}Z.
$$
With respect to this basis the induced metric and its inverse are
$$
g=\left(\begin{array}{cc}
\lambda+(f')^2 & -\frac{\lambda}{2}f^2\\
-\frac{\lambda}{2}f^2 & f^2+\frac{\lambda}{4}f^4
\end{array}\right),\quad g^{-1}=f^{-2}D^{-1}\left(\begin{array}{cc}
 f^2+\frac{\lambda}{4}f^4& \frac{\lambda}{2}f^2\\
\frac{\lambda}{2}f^2 & \lambda+(f')^2
\end{array}\right),
$$
where $D=\lambda+\frac{\lambda}{4}f^2(f')^2+(f')^2$. The unit normal vector is 
$$
\nuv=-\frac{\sql}{\sqrt{D}}\left(\frac{ff'}{2}\sss+\ccc\right)X+\frac{\sql}{\sqrt{D}}\left(\frac{ff'}{2}\ccc-\sss\right)Y+\frac{f'}{\sqrt{\lambda D}}Z.
$$
By Lemma \ref{2ff} we get:
\begin{eqnarray*}
\bar{\nabla}_{V_1}V_1 & = & \left(f''\ccc+\sql f'\sss\right)X+\left(f''\sss-\sql f'\ccc\right)Y;\\
\bar{\nabla}_{V_2}V_2 & = & -f\left(1+\frac{\sql}{2}f^2\right)\left(\ccc X+\sss Y\right);\\
\bar{\nabla}_{V_1}V_2 & = & \left(-f'\sss+\frac{\sql}{4}\left(2f\ccc-f^2f'\sss\right)\right)X\\
&&+\left(f'\ccc+\frac{\sql}{4}\left(2f\sss+f^2f'\ccc\right)\right)Y-\frac{ff'}{2}Z.
\end{eqnarray*}
After some standard computations the expression of the second fundamental form  and of the mean curvature are
$$
\begin{array}{l}
A=\frac{\sql}{\sqrt D}\left(\begin{array}{cc}
-f''-\frac{\sql}{2}f(f')^2 & \frac{\sql}{8}f\left(f^2(f')^2-4\right)\\
\frac{\sql}{8}f\left(f^2(f')^2-4\right) & f\left(1+\frac{\sql}{2}f^2\right)
\end{array}\right),\\
H =  tr(Ag^{-1})=\frac{\sql}{f D^{\frac 32}}\left(-ff''\left(1+\frac{\lambda}{4}f^2\right)+\lambda+(f')^2\right).
\end{array}
$$
The vertical translating condition \eqref{vert trasl} becomes
$$
\frac{\sql}{f D^{\frac 32}}\left(-ff''\left(1+\frac{\lambda}{4}f^2\right)+\lambda+(f')^2\right)=\frac{f'}{\sqrt{D}}.
$$
Following the strategy of \cite{CSS} we fix any $f_0>0$ and we look for the solution of the following Cauchy problem
\begin{equation}\label{ode6}
\left\{\begin{array}{rcl}
f''&=&\frac{4}{\sql f\left(4+{\lambda}f^2\right)}\left(\lambda^{\frac 32}+\sql (f')^2-f(f')^3-\lambda ff'-\frac{\lambda}{4}f^3(f')^3\right).\\
f(0)&=&f_0,\quad f'(0)\ =\ 0.
\end{array}\right.
\end{equation}
Such $f$ is defined at least in a neighborhood of $0$. Moreover, by \eqref{ode6} we have
$$
f''(0)=\frac{4\lambda}{f_0\left(4+\lambda f^2_0\right)}>0,
$$ 
then we can find an $\varepsilon>0$ such that $f$ is concave in the interval $(-2\varepsilon,2\varepsilon)$ and $f'>0$ (resp. $f'<0$) in $\left(0,2\varepsilon\right)$ (resp. in $\left(-2\varepsilon,0\right)$). Define the following four constants: $c_0^{\pm}=f(\pm\varepsilon)$ and $c_1^{\pm}=f'(\pm\varepsilon)$. Let $\varphi^{\pm}$ be the solutions of \eqref{ode rotazione} with initial conditions
$$
\varphi^{\pm}(c_0^{\pm})=\pm\varepsilon,\quad (\varphi^{\pm})'(c_0^{\pm})=\frac{1}{c_1^{\pm}}.
$$
From what proved so far, $\varphi^{\pm}$ is defined for every $r$ bigger than $c_0^{\pm}$, hence gluing the graphs of $f$, $\varphi^+$ and $\varphi^-$ we have a complete unbounded curve $\mathcal R=\mathcal R(f_0,\lambda)$. Moreover this curve is embedded because of the uniqueness of the solution of the ODE and $f_0$ represent the distance of this curve from the $z$-axis. The rotation of this curve along the $z$-axis produces a complete vertical translators. On the other hand every vertical translator invariant by horizontal rotations can be found in this way up to isometries of $\hhh_3$.
\cvd

\begin{figure}[H]
\centering
\includegraphics[width=0.6\textwidth]{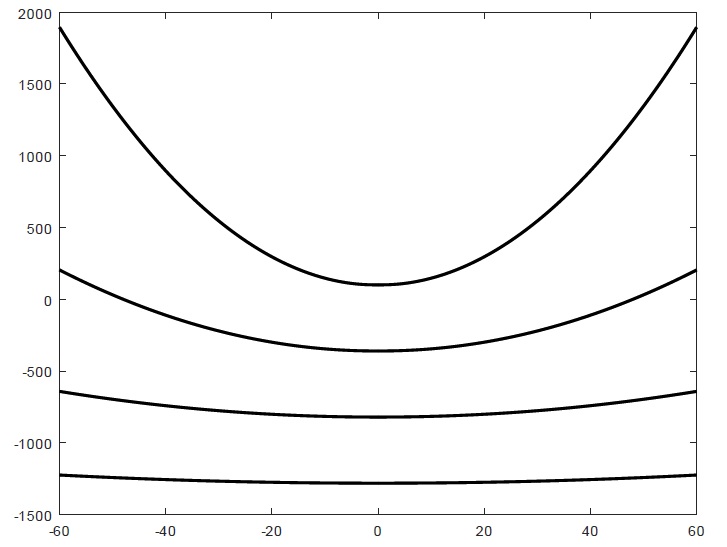}
\caption{Some examples of graphs generating the bowl solution in $\hhh_3$: from top to bottom $\lambda=1$, $\lambda=10$, $\lambda=100$ and $\lambda=1000$. The curves have been translated vertically in order to avoid overlapping.}
\end{figure}

\begin{figure}[H]
\centering
\includegraphics[width=0.6\textwidth]{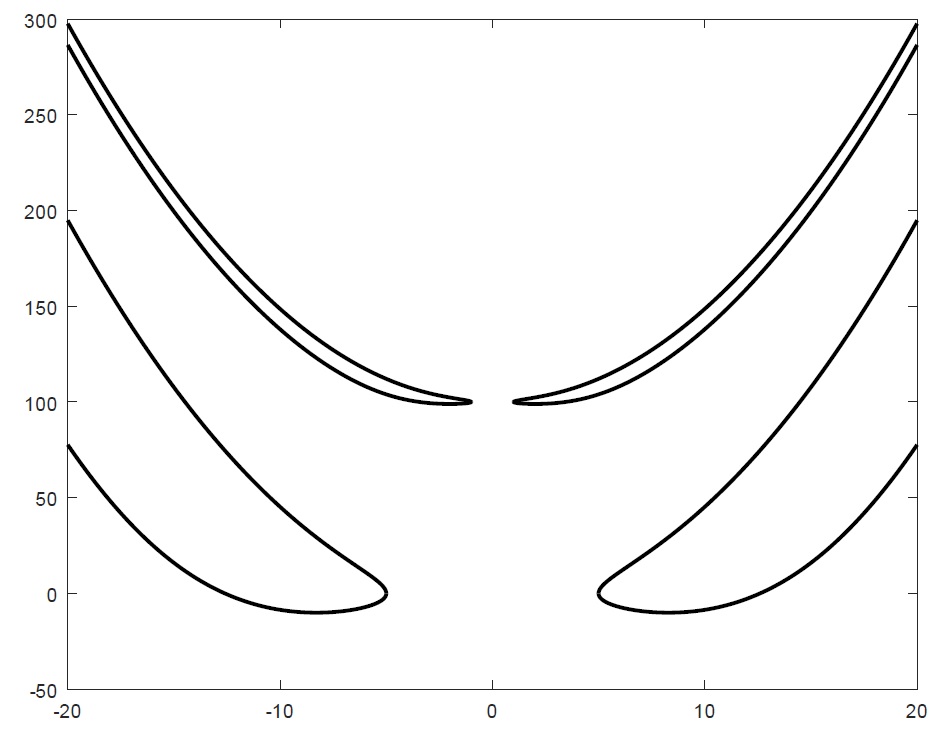}
\caption{Some examples of curves generating translating catenoids: $\lambda=1$ in every cases, $f_0=1$ above and $f_0=5$ below. The curves have been translated vertically in order to avoid overlapping.}
\end{figure}

Now we want to understand the asymptotic behavior of each arm of the translators described in Theorem \ref{teor rot}.

\begin{Lemma}\label{rotazione stime}
The function $\varphi$ solution of equation \eqref{ode rotazione} is defined for any $r$ sufficiently big and it has the following asymptotic expansion as $r$ goes to infinity:
$$
\varphi(r)\approx \frac{r^2}{2\sqrt{\lambda}}+\left\{\begin{array}{ccc}
-\frac{4}{\sqrt{\lambda}(4-\lambda)}\log(r) & \text{if} & \lambda<4;\\
-\log^2(r) & \text{if} & \lambda=4;\\
C_0 r^{1-\frac{4}{\lambda}} & \text{if} & \lambda>4,
\end{array}\right.
$$
for some $C_0$ depending on the initial datum.
\end{Lemma}
\proof
Although we will be forced to split into cases, until it is explicitly said $\lambda$ is any positive number. Let $\psi=\varphi'$, by \eqref{ode rotazione} we have
\begin{equation}\label{ode1}
\psi'=\frac{1}{\sqrt{\lambda}}+\frac{4\psi}{r(\lambda r^2+4)}\left(\sqrt{\lambda} r \psi-1-\lambda\psi^2\right).
\end{equation}
$\psi'$ vanishes on the set:
$$
\mathcal C^{\lambda}=\left\{r(\lambda r^2+4)+4\sqrt{\lambda}\psi\left(\sqrt{\lambda}r\psi-1-\lambda\psi^2\right)=0\right\}.
$$
$\mathcal C^{\lambda}$ can be seen as an algebraic curve in the plane of coordinates $(r,\psi)$ depending on a parameter $\lambda$. When $r$ goes to infinity, this curve is asymptotic to a straight line of equation
$$
r=\bar R\psi,
$$
where $\bar R=\bar R(\lambda)$ is the unique real solution of the equation $x^3+4x=4\sqrt{\lambda}$. It can be checked that $\bar R$ is positive and diverges as $\lambda$ goes to infinity as $\lambda^{\frac 16}$. This implies that $\mathcal C^{\lambda}$ is approaching the horizontal axis when $\lambda$ diverges. Moreover $\psi'$ is negative above $\mathcal C^{\lambda}$ and positive below. 
Using the implicit function theorem we can see that the curve $\mathcal C^{\lambda}$ is the graph of a strictly increasing function, then it behaves as a barrier for $\psi$: once $\psi$ enters in the region below $\mathcal C^{\lambda}$, it will stay there forever. Therefore we can deduce that $\psi$ is defined for any $r\geq r_0$, it grows at most linearly with a slope that depends on $\lambda$ and there exists a $\bar r> r_0$ such that for any $r\geq \bar r$ we have that $\psi'>0$. Moreover we have that $\psi$ becomes positive at least for large $r$. In fact, suppose that there exists some $\varepsilon>0$ such that $\psi<-\varepsilon$ for any $r\geq\bar r$, then we have that, up to increase $\bar r$, for any $r\geq\bar r$:
$$
\psi'>\frac{1}{\sqrt{\lambda}}-\frac{1}{r^2}.
$$
Integrating we can see that $\psi$ has to become positive for a finite $r\geq 0$, having a contradiction. Next we show that $\psi$ grows like $\frac{r}{\sqrt{\lambda}}$: let us define the function $\gamma=\frac{\psi}{r}$, the goal is to prove that $\gamma$ converges to $\frac{1}{\sqrt{\lambda}}$. By \eqref{ode1} we get:
\begin{eqnarray}
\nonumber\gamma'& = & \frac{\psi'}{r}-\frac{\psi}{r^2}=\frac 1r\left[\frac{1}{\sql}+\frac{4\gamma}{\lambda r^2+4}\left(\sql r^2\gamma-1-\lambda r^2\gamma^2\right)-\gamma\right]\\
&=&\frac 1r\left(4\gamma^2+1\right)\left(\frac{1}{\sql}-\gamma\right)+\frac{4\gamma}{r\left(\lambda r^2+4\right)}\left(4\gamma^2-\frac{4}{\sql}\gamma-1\right).\label{ode2}
\end{eqnarray}
We know that $\gamma$ is positive and bounded. Moreover, by \eqref{ode2} we can see that if $\gamma=\frac{1}{\sql}$, then $\gamma'<0$, therefore, increasing $\bar r$ if necessary, we can suppose that either $\gamma>\frac{1}{\sql}$ or $\gamma<\frac{1}{\sql}$ for any $r\geq\bar r$. It follows that there are positive constants $b_1$ and $b_2$ such that 
$$
-\frac{b_1}{r}\left(\gamma-\frac{1}{\sql}\right)-\frac{1}{r^2}<\gamma'<-\frac{b_2}{r}\left(\gamma-\frac{1}{\sql}\right)+\frac{1}{r^2}
$$
holds for any $r\geq\bar r$ (once again we can increase $\bar r$ if necessary). Applying twice Lemma \ref{ode speciale} we can estimate $\gamma$ from both sides with two functions that converge to the constant $\lambda^{-\frac 12}$, in particular we deduce that also $\gamma$ converges to $\lambda^{-\frac 12}$ as claimed.

Now we consider the function $q=\psi-\frac{r}{\sql}$. We claim that $q$ converges to $0$.  First we can compute its evolution equation using \eqref{ode1}. 
\begin{eqnarray}
 q' & = &-\frac {4}{r(\lambda r^2+4)}\left(q+\frac{r}{\sql}\right)\left(\lambda q^2+\sql rq+1\right). \label{ode3}
\end{eqnarray}
Moreover it follows from what has been seen so far that $q$ is sublinear. We can improve this first estimate saying that $q$ is bounded. In fact if $q\geq 0$, then clearly $q'<0$, hence $q$ is bounded from above, and if $q$ becomes negative, then it stays negative forever. Suppose now that $q<0$ and $q'<0$ for every $r\geq\bar r$. In view of the sublinearity we can see that for every $\varepsilon>0$ 
\begin{equation}\label{xxx1}
-\varepsilon r<q<-\varepsilon
\end{equation}
holds for every $r$. By \eqref{ode3} and the sublinearity we have that
\begin{eqnarray}
q'<0&\Leftrightarrow&\lambda q^2+\sql rq+1>0\label{xxx2}
\end{eqnarray}
but \eqref{xxx1} and \eqref{xxx2} are in contradiction if $r$ is too big. Therefore there is a positive $\bar r$ such that for every $r\geq\bar r$ we have $qq'<0$. In particular $q$ is bounded from below too and it becomes monotone after $\bar r$. Finally we can prove that $q$ converges to zero. Suppose now that there exists a positive constant $\varepsilon$ such that $q\geq\varepsilon$ for every $r\geq\bar r$. In this case, by \eqref{ode3}, we get that
\begin{eqnarray*}
q'&\leq&-\frac {4}{r(\lambda r^2+4)}\left(\varepsilon+\frac{r}{\sql}\right)\left(\lambda \varepsilon^2+\sql r\varepsilon+1\right) \leq-\frac cr,
\end{eqnarray*}
for some positive constant $c$ if $r\geq\bar r$. Integrating we have that $q$ becomes negative for finite $r$, having a contradiction. On the other hand, if we suppose that there exists a positive constant $\varepsilon$ such that $q\leq-\varepsilon$ for every $r\geq\bar r$, from the fact that $q$ is bounded, we have that, up to increase $\bar r$, we can estimate \eqref{ode3} as follows:
\begin{eqnarray*}
q'&\geq&-\frac {4}{r(\lambda r^2+4)}\left(-\varepsilon+\frac{r}{\sql}\right)\left(\lambda \varepsilon^2-\sql r\varepsilon+1\right)\geq\frac cr,
\end{eqnarray*}
for some positive constant $c$ if $r$ is big enough. Arguing as before we get a contraddiction. This shows that $q$ converges to $0$ as claimed.

Now we want to conclude the proof of this Lemma studying how fast $q$ converges to $0$. In this case the value of $\lambda$ has a crucial role. Fix a positive constant $k$ and let us define the following function: $\eta=q r^k$. By \eqref{ode3} we can compute:
\begin{eqnarray}
\eta' &=& \frac{1}{r^{2k+1}(\lambda r^2+4)}\left((k\lambda-4)\eta r^{2k+2}+4(k-1)\eta r^{2k}\phantom{\frac{4}{\sql}}\right.\label{ode4}\\
\nonumber &&\phantom{aaaaaaaaaaaaaaa}\left.-4\lambda\eta^3-8\sql r^{k+1}\eta^2-\frac{4}{\sql}r^{3k+1}\right).
\end{eqnarray}
Now we have to split the proof into three cases according if $\lambda$ is smaller, equal or bigger than $4$.
\begin{itemize}
\item[i)] Suppose $\lambda<4$. Let in this case $k=1$. We claim that $\eta$ converges to $-\frac{4}{\sql(4-\lambda)}$. Equation \eqref{ode4} can be written in the following way:
\begin{eqnarray}
\nonumber\eta' & = & -\frac{4-\lambda}{\lambda r}\left(\eta+\frac{4}{\sql(4-\lambda)}\right)+\frac{4(4-\lambda)}{\lambda r(\lambda r^2+4)}\left(\eta+\frac{4}{\sql(4-\lambda)}\right)\\
&&-\frac{4\lambda\eta^3}{r^3(\lambda r^2+4)}-\frac{8\sql\eta^2}{r(\lambda r^2+4)}\label{rotBB}
\end{eqnarray}

Since we proved that $q$ converges to $0$, then the function $\eta$ is sublinear. By \eqref{rotBB} and the sublinearity, for every $\varepsilon>0$ there are two positive constant $C_1$ and $C_2$ such that
$$
-\frac{4-\lambda}{\lambda r}\left(\eta+\frac{4}{\sql(4-\lambda)}+\varepsilon\right)+\frac{C_2}{r^2}<\eta'<-\frac{4-\lambda}{\lambda r}\left(\eta+\frac{4}{\sql(4-\lambda)}\right)+\frac{C_1}{r^2}.
$$
Applying twice Lemma \ref{ode speciale} and letting $\varepsilon$ go to zero, we can prove our claim.

\item[ii)] Suppose now that $\lambda>4$. This time take $k=\frac{4}{\lambda}$. Since we already know that $q$ converges to $0$, we have that for any positive constant $c$ there is a $\bar r$ big enough such that $|\eta|<cr^k$ for any $r\geq\bar r$. Using this information and equation \eqref{ode4} we have that there is positive constant $b$ such that
$$
|\eta'|<br^{k-2}
$$
holds for every $r$ big enough. Since $k<1$, integrating we can see that $\eta$ is bounded. By \eqref{ode4} we can see that function $\eta'$ vanishes on the curve of equation
$$
\lambda\eta^3+2\sql r^{k+1}\eta^2+(1-k)\eta r^{2k}+\frac{r^{3k+1}}{\sql}=0.
$$
It can be thought as a polynomial of degree three in the variable $\eta$. Can be proved that this polynomial has an unique real root if $r$ is big enough and that it is negative.  Moreover $\eta'<0$ above the curve and $\eta'>0$ below. Hence there is a $\bar r$ such that $\eta$ is increasing for $r\geq\bar r$. Since $\eta$ is monotone and bounded, it follows that it converges to some constant $C_0$ wich depends on the initial value of $\eta$.
\item[iii)] Finally we consider the case $\lambda=4$. Let $k=1$, equation \eqref{ode4} becomes:
\begin{equation}\label{ode4iii}
\eta' = -\frac{(8\eta^3+8r^2\eta^2+r^4)}{2r^3(r^2+1)}.
\end{equation}
Since $q$ converges to $0$, $\eta$ is sublinear in this case too. Unfortunally this time $\eta$ is unbounded. In fact, from \eqref{ode4iii} and the sublinearity of $\eta$, we can find two positive constants $c_1,c_2$ such that  the following estimate holds for every $r\geq\bar r$:
$$
\frac{c_1}{r^2}-\frac{1}{r}\leq\eta'\leq\frac{c_2}{r^2}-\frac{1}{4r}.
$$
Integrating we have that there are some constants $c_3,c_4$ such that
\begin{equation}\label{yyy1}
c_3-\frac{c_1}{r}-\log(r)\leq\eta\leq c_4-\frac{c_2}{r}-\frac{\log(r)}{4}.
\end{equation}
In particular it follows that $\eta$ diverges as $r$ goes to infinity. Let us define the function $\zeta=\frac{\eta}{\log(r)}$. By \eqref{ode4}, its evolution equation is:
\begin{eqnarray}
\zeta' &=&-\frac{2\zeta+1}{2r\log(r)}+\frac{1}{2\log(r)r(r^2+1)}-\frac{4\log^2(r)\zeta^3}{r^3(r^2+1)}-\frac{4\log(r)\zeta^2}{r(r^2+1)}.\label{ode5}
\end{eqnarray}
Estimate \eqref{yyy1} says that $\zeta$ is bounded and negative for large values of $r$, then up to increase $\bar r$, we get:
$$
-\frac{1}{r^2}-\frac{2\zeta+1}{2r\log(r)}\leq\zeta'\leq-\frac{2\zeta+1}{2r\log(r)}+\frac{1}{r^2}.
$$
It follows that there are two constants $c_5$ and $c_6$ such that
$$
-\frac 12+\frac 1r+\frac{1}{r\log(r)}+\frac{c_5}{\log(r)}\leq\zeta\leq\frac{c_6}{\log(r)}-\frac 1r-\frac{1}{r\log(r)}-\frac 12.
$$
In particular we have that $\zeta$ converges to $-\frac 12$. 
\end{itemize}
In thes way we found the asymptotic behavior of $\psi=\varphi'$. Integrating we get the thesis.\cvd


Now we want to investigate the convergence when $\lambda$ diverges. We can see that the equation \eqref{ode rotazione} converges to
$$
\varphi''=-4\frac{(\varphi')^3}{r^3}.
$$
The constant functions are trivial solutions of this equations. On the other hand \eqref{ode6} converges to
$$
\left\{\begin{array}{l}
f''  =  \frac{4}{f^3},\\
f(0)=f_0,\ f'(0)=0.
\end{array}\right.
$$
It can be easily checked that the explicit solutions are $\tfz=\frac{\sqrt{4z^2+f_0^4}}{f_0}$. We wish to say that the rotationally invariant vertical translators converge to either a horizontal plane or the surface of equation $\sqrt{x^2+y^2}=\frac{\sqrt{4z^2+f_0^4}}{f_0}$. In fact part $2)$ and $3)$ of Theorem \ref{main 2} hold as a direct corollary of the following two results.


\begin{figure}[H]
\centering
\includegraphics[width=0.6\textwidth]{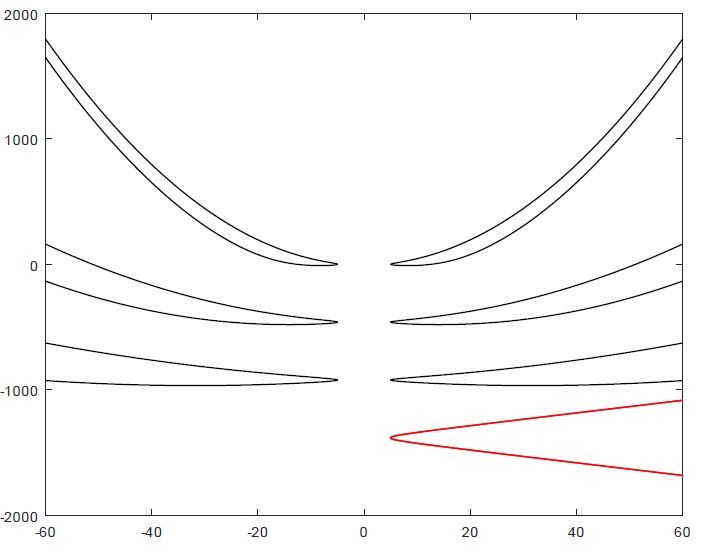}
\caption{Some examples curve generating a translating catenoids: $f_0=5$ in every cases, from top to above $\lambda=1$, $\lambda=10$ and $\lambda=100$. In red we have the graph of limit curve $\sqrt{x^2+y^2}=\frac{\sqrt{4z^2+5^4}}{5}$. The curves have been translated vertically in order to avoid overlapping.}
\end{figure}

\begin{Proposition}
Let $\psi$ be unique entire solution of \eqref{ode1}, then there is a constant $C_0$ such that 
$$
\left|\psi\right|\leq C_0\lambda^{-\frac 16},
$$
and for every $m\in\mathbb N$ there is a positive constant $C_m$ such that 
$$
\left|\psi^{(m)}\right|\leq C_m\lambda^{-\frac 12}.
$$
\end{Proposition}
\proof The case $m=0$ follows from what saw in the proof of Lemma \ref{rotazione stime}: $\psi\geq 0$ for every $r$, its graphs stays below the curve $\mathcal C^{\lambda}$ and this curve is asymptotic to the straight line $r=c\lambda^{\frac 16}z$ for some constant $c$.
The case $m=1$ can be easly derived by the evolution equation for $\psi$ \eqref{ode1}, the estimate for $m=0$ and the fact that $\lim_{r\rightarrow 0}\frac{\psi}{r}=\frac{1}{2\sql}$.

Now fix an $m\geq 1$. We define for brevity of notation $\zeta=\zeta(r)=4+\lambda r^2$. Since for every $j$
$$
\left(\frac 1r\right)^{(j)} = \frac{(-1)^j j!}{r^{j+1}},
$$
by Lemma \ref{derivate} and \eqref{ode1} we can prove that
\begin{equation}\label{NUM}
\psi^{(m+1)}  =  4\sum_{j=0}^m \sum_{k=0}^{m-j}\left[{{m}\choose{j}}{ {m-j}\choose{k} }\frac{(-1)^jj!}{r^{j+1}}\left(\frac{1}{\zeta}\right)^{(m-j)}\left(\sql r\psi^2-\psi-\lambda\psi^3\right)^{(k)}\right]
\end{equation}
By Lemma \ref{derivate} and the explicit expression of $\zeta$ we can see that for every $j$
$$
\left|\left(\frac{1}{\zeta}\right)^{(j)}\right|\leq c_j\lambda^{-1},
$$
uniformly in $r$ for some constant $c_j$. Moreover, using again Lemma \ref{derivate}, we have
\begin{equation}\label{fff}
\begin{array}{l}
\left(\sql r\psi^2-\psi-\lambda\psi^3\right)^{(k)} = \\
\phantom{aaaa}= -\psi^{(k)}+\sql r\sum_{l=0}^k{k\choose l}\psi^{(l)}\psi^{(k-l)}+k\sql\sum_{l=0}^{k-1}\psi^{(l)}\psi^{(k-l-1)}\\
\phantom{aaaa=} -\lambda\sum_{l=0}^k\sum_{q=0}^{k-l}{k\choose l}{{k-l}\choose q}\psi^{(l)}\psi^{(q)}\psi^{(k-l-q)}.
\end{array}
\end{equation}
Therefore, in view of the estimate for $\psi$ given above and supposing that the thesis holds for every number smaller or equal to $m$, we get that for every $j$ and $k$:
$$
\left|\frac{\left(\sql r\psi^2-\psi-\lambda\psi^3\right)^{(k)}}{r^{j+1}}\right|\leq c_{j,k}\lambda^{\frac 12},
$$
for some positive constant $c_{j,k}$. The worst case is the addend in \eqref{fff} with no derivatives. The apparent singularity in $r=0$ can be eliminated noting that for every $m\in\mathbb N$ there is a constant $\tilde c_m$ that can be computed recursively from \eqref{NUM} such that, when $r\rightarrow 0$,
$$
\psi^{(2m)}\rightarrow 0,\quad\psi^{(2m+1)}\rightarrow \frac{\tilde c_m}{\sql}.
$$
\cvd

\begin{Proposition}
Fix $f_0>0$, then for every compact $K\subset\mathbb R$ there is a $\lambda_0>0$ such that $f$ - solution of \eqref{ode6} - is defined on $K$ for every $\lambda\geq\lambda_0$. Moreover for every $m\in\mathbb N$ there is a positive constant $C$ depending only on $K$, $m$ and $f_0$ such that
$$
\left|f^{(m)}-\tf^{(m)}\right|(z)\leq C\lambda^{-\frac 12}
$$
for every $z\in K$ and $\lambda\geq\lambda_0$.
\end{Proposition}
\proof
Fix any compact $K$. We start proving that all the derivatives of $f$ and $\tf$ are bounded on $K$ uniformly in $\lambda$. After that we will improve these estimates showing that all the derivatives of the function $\psi=f-\tf$ converge to zero as $\lambda^{-\frac 12}$.
It is trivial to compute that $\tf$ and $\tf'$ are bounded in $K$. Since $\tf''=4\tf^{-1}$ and $\tf\geq f_0$, it is easy to prove by induction that all the derivatives of $\tf$ are bounded in $K$ too. The proof for $f$ is more involved because $f$ depends on $\lambda$. From the construction of the curve $\mathcal R(f_0,\lambda)$ in the proof of Theorem \ref{teor rot}, we can see that $f$ has a global minimum in $z=0$. Therefore by \eqref{ode6} we can estimate
\begin{equation}\label{zpn}
f''\ \leq\ \left\{\begin{array}{rcl}
\frac{4}{f_0^3}\left(1+\frac{1}{\sql}(f')^2\right) & \text{if} & z\geq 0,\\
\frac{8}{f_0^3}-\frac{2}{\sql}(f')^3 &\text{if}& z\leq 0\ \text{and}\  \lambda\  \text{big enough}.
\end{array}\right.
\end{equation}

From the first estimate in \eqref{zpn} and what we saw so far we can deduce that if $z\geq 0$ we have
\begin{eqnarray*}
0\ \leq\  f'(z)&\leq&\sql\tan\left(\frac{4z}{f_0^3\sql}\right),\\
f_0\ \leq\  f(z) & \leq & f_0-\frac{f_0^3\lambda}{4}\log\left(\cos\left(\frac{4z}{f_0^3\sql}\right)\right).
\end{eqnarray*}
From these estimates we can see that $f(z)$ is defined for any fixed $z>0$ if $\lambda$ is big enough and that $f$ and $f'$ are bounded in $K\cap [0,+\infty)$ with a constant not depending on $\lambda$. If $z\leq 0$, by the second estimate in \eqref{zpn} we have that $f'(z)\geq\varphi(z)$ where $\varphi$ is solution of
\begin{equation}\label{ode66}
\varphi'=a^3-b^3\varphi^3,\quad \varphi(0)=0
\end{equation}
with, for brevity of notation, $a={2}{f_0}^{-1}$ and $b=2^{\frac 13}\lambda^{-\frac 16}$. Standard computations show that $\varphi$ solves the following implicit equation
\begin{equation}\label{implicita}
z=\frac{1}{3ab^2}\log\left(\frac{\sqrt{b^2\varphi^2+ab\varphi+a^2}}{a-b\varphi}\right)+\frac{1}{\sqrt{3}a^2b}\left(\arctan\left(\frac{a+2b\varphi}{a\sqrt{3}}\right)-\arctan\left(\frac{1}{\sqrt{3}}\right)\right).
\end{equation}
Therefore $\varphi$ blows up at $z=z_{\lambda}:=-\frac{2\sqrt{3}\pi}{9a^2b}=-2^{-\frac43}3^{-\frac 32}\pi f_0^2\lambda^{\frac 16}$, and $z_{\lambda}$ diverges with $\lambda$. Hence $f$ is defined at any fixed negative $z$ if $\lambda$ is big enough. Let $z_0$ be the infimum of $K$ and suppose that it is negative. By \eqref{ode66} we get that $\varphi$ is increasing, therefore it is sufficient to estimate $\varphi_0=\varphi(z_0)$ in order to estimate $\varphi$ and hence $f'$ on $K\cap(-\infty,0]$. We know that $\varphi_0<0$, then by \eqref{implicita} we have
$$
\frac{1}{3ab^2}\log\left(\frac{\sqrt{b^2\varphi_0^2+ab\varphi_0+a^2}}{a-b\varphi_0}\right)>z_0.
$$
Since $\varphi(0)=0$, standard computations lead to say that the previous inequality holds if and only if
$$
\varphi_0>-\frac{12z_0}{f_0^3}\frac{\sqrt{3\mu(4-\mu)}-\mu-2}{\log(\mu)(\mu-1)},
$$
where $\mu=\mu(\lambda,z_0,f_0)=e^{-3a^2cz_0}=e^{-2^{\frac 73}3f_0^{-2}\lambda^{-\frac 16}z_0}$. Since for every $z_0$ and $f_0$ 
$$
\lim_{\lambda\rightarrow\infty}\frac{\sqrt{3\mu(4-\mu)}-\mu-2}{\log(\mu)(\mu-1)}=-\frac 23,
$$
we can bound $\varphi_0$ in $K\cap(-\infty,0]$ uniformly in $\lambda$ if $\lambda$ is sufficiently large. It follows for example that if $\lambda$ is big enough for every $z\in K\cap(-\infty,0]$
\begin{eqnarray*}
\frac{16z_0}{f_0^3}\ \leq\ \varphi_0\ \leq\ \varphi(z)\ \leq&f'(z)& \leq\  0,
\end{eqnarray*}
therefore
\begin{eqnarray*}
f_0\ \leq&f(z)& \leq f_0+\frac{16z_0^2}{f_0^3}.
\end{eqnarray*}
Summarizing what has been seen so far, on any compact $K$ we are able to bound $f$ and $f'$ uniformly in $\lambda$. By \eqref{ode6}, $f''$ is bounded too. Applying Lemma \ref{derivate} to \eqref{ode6} we can prove by induction that for every $K$ and $m\in\mathbb N$ there is a constant $c>0$ depending on $K$, $m$ and $f_0$ such that $\left|f^{(m)}\right|<c$.
Let us define the following function $\psi=f-\tf$. By the evolution equations of $f$ and $\tf$ we have that
\begin{equation}\label{odepsi}
\left\{\begin{array}{rcl}
\psi'' &=& -4\psi\left(f^2+f\tf+\tf^2\right)f^{-3}\tf^{-3}+\frac{4\left(\sql f^2 (f')^2-f^3(f')^3-\lambda f^3f'-\frac{\lambda}{4}f^5(f')^3-4\sql\right)}{\sql f^3\left(4+{\lambda}f^2\right)}\\
\psi(0)&=&\psi'(0)\ =\ 0
\end{array}\right.
\end{equation}
From \eqref{odepsi} and the bounds for $f$, $f'$ and $\tf$ saw before we have that $\psi$ is a solution of the following differential inequality if $\lambda$ is big enough
\begin{equation}\label{diffineq}
\left|\psi''\right| \leq c\left(\left|\psi\right|+\lambda^{-\frac 12}\right),\quad \psi(0)\ =\ \psi'(0)\ =\ 0
\end{equation}
for some constant $c$ depending only on $K$ and $f_0$. By Theorem 21.1 of \cite{Sz} we have that
\begin{equation}
\left|\psi(z)\right|\leq \frac{1}{\sql}\left(\cosh(\sqrt{c}z)-1\right),\quad\left|\psi'(z)\right|\leq \frac{\sqrt{c}}{\sql}\sinh(\sqrt{c}|z|).
\end{equation}
Therefore the thesis follows for $m=0$ and $m=1$. Using this information and \eqref{diffineq} we can prove the thesis for $m=2$ too. Let $m$ be any integer greater than $1$ and suppose that the thesis holds for every $i\leq m$. By \eqref{odepsi} and Lemma \ref{derivate} we get
\begin{eqnarray*}
\psi^{(m+2)} & = & -4\sum_{i=0}^m{m \choose i}\psi^{(m-i)}\left[(f^2+f\tf+\tf^2)f^{-3}\tf^{-3}\right]^{(i)}\\
&&+4\left[\frac{\left(\sql f^2 (f')^2-f^3(f')^3-\lambda f^3f'-\frac{\lambda}{4}f^5(f')^3-4\sql\right)}{\sql f^3\left(4+{\lambda}f^2\right)}\right]^{(m)}.
\end{eqnarray*}
The terms in square brackets can be expanded using Lemma \ref{derivate} again. Since we proved that all the derivatives of $f$ and $\tf$ are bounded we can see that the first one is bounded, while the second one decays as $\lambda^{-\frac 12}$. The thesis follows applying the inductive hypothesis on the terms $\psi^{(m-i)}$ for every $i\leq m$.
\cvd

\section{Helicoidal surfaces}
In this section we finish the proof of Theorem \ref{main 1} examining the helicoidal vertical translators. We start with a helicoidal surface $\mm$ with vertical axis. There exists a planar curve $\gamma=(\gamma_1,\gamma_2)$ such that $\mm$ can be parametrized in the following way:
\begin{eqnarray*}
\mm(v_1,v_2)&=&\left(\gamma_1(v_1)\cos(v_2)-\gamma_2(v_1)\sin(v_2),\gamma_1(v_1)\sin(v_2)+\gamma_2(v_1)\cos(v_2),cv_2\right)\\
&=&\left(e^{iv_2}\gamma(v_1),cv_2\right),
\end{eqnarray*}
where $c\in\mathbb R$ is the pitch of the helicoidal motion and, in the last equality, we use complex coordinates on $\hhh_3\equiv \mathbb C\times\mathbb R$. Suppose that $\gamma$ is parametrized by arc length, we introduce the notations of \cite{Ha} useful to simplify some computations.  Let $T=\gamma'=\frac{\partial\gamma}{\partial v_1}$ be the unit tangent vector of the curve, hence $N=iT=(-\gamma_2',\gamma_1')$ is the unit normal vector field. Let us denote with $\left\langle \cdot,\cdot\right\rangle$ the usual Euclidean scalar product on the plane, then we define the following functions: 
\begin{eqnarray*}
\tau&=&\left\langle\gamma,T \right\rangle =\gamma_1\gamma_1'+\gamma_2\gamma_2';\\
\nu &=& \left\langle \gamma, N\right\rangle=\gamma_1'\gamma_2-\gamma_1\gamma_2';\\
r^2&=& \left\langle \gamma,\gamma\right\rangle=\tau^2+\nu^2;\\
k &=& \left\langle \gamma'',N\right\rangle=-\gamma_1''\gamma_2'+\gamma_1'\gamma_2''.
\end{eqnarray*}
Obviusly $k$ is the curvature of $\gamma$ and $\gamma=\tau T+\nu N=(\tau+i\nu)T$. Deriving with respect to the arc length we have:
\begin{equation}\label{tau-nu}
\left\{\begin{array}{rcl}
\tau' & = & 1+k\nu,\\
\nu' & = & -k\tau.
\end{array}\right.
\end{equation}

\begin{Lemma}\label{mc-helicoidal}
For any given $c\in\mathbb R$, $\lambda>0$ and curve $\gamma$, the mean curvature of $\mm$ in $\left(\hhh_3,\bar g_{\lambda}\right)$ is 
$$
H=\frac{2\sql \left[-ck\left(4r^2+{\lambda}(2c-r^2)^2\right)+{\lambda c}\nu\left(4c-2\tau^2-\nu^2\right)-4c\nu\right]}{\left(4\tau^2+{\lambda}\left(4c^2-4c\tau^2+\tau^4+\tau^2\nu^2\right)\right)^{\frac 32}}.
$$
Furthermore $\mm$ is a vertical translator with unit speed if and only if the curvature $k$ is
\begin{equation}\label{eq k}
k=\frac{{\lambda^{\frac 32} c}\nu\left(4c-2\tau^2-\nu^2\right)-4\sql c\nu-\tau\left(4\tau^2+{\lambda}\left(4c^2-4c\tau^2+\tau^4+\tau^2\nu^2\right)\right)}{\sql c\left(4r^2+{\lambda}\left(2c-r^2\right)^2\right)}.
\end{equation}
\end{Lemma}
\proof
With the notations of Lemma \ref{2ff} we have:
$$
\begin{array}{rclcrcl}
a_1 & = & \gamma_1'(v_1)\cos(v_2)-\gamma_2'(v_1)\sin(v_2)=:x_1, & \quad & a_2 & = & -y,\\
b_1 & = & \gamma_1'(v_1)\sin(v_2)+\gamma_2'(v_1)\cos(v_2)=:y_1, & \quad & b_2 & = & x,\\
c_1 & = & \frac{\nu}{2}, & \quad & c_2 & = & c-\frac{r^2}{2}.
\end{array}
$$
It follows that the tangent space of $\mm$ is generated by
$$
V_1=x_1 X+y_1 Y+\frac{\nu}{2}Z,\quad V_2=-yX+xY+\frac{1}{2}(2c-r^2)Z=F_4+cF_3.
$$
Moreover it is easy to check:
$$
\begin{array}{rclcrcl}
\frac{\partial a_1}{\partial v_1} & = & \gamma_1''(v_1)\cos(v_2)-\gamma_2''(v_1)\sin(v_2)=:x_{11},& \quad & \frac{\partial a_1}{\partial v_2} & = & -y_1,\\
\frac{\partial b_1}{\partial v_1} & = & \gamma_1''(v_1)\sin(v_2)+\gamma_2''(v_1)\cos(v_2)=:y_{11},& \quad & \frac{\partial b_1}{\partial v_2} & = & x_1,\\
\frac{\partial c_1}{\partial v_1} & = & -\frac{k\tau}{2},& \quad & \frac{\partial c_1}{\partial v_2} & = & 0,\\
\frac{\partial a_2}{\partial v_1} & = & -y_1,& \quad & \frac{\partial a_2}{\partial v_2} & = & -x,\\
\frac{\partial b_2}{\partial v_1} & = & x_1,& \quad & \frac{\partial b_2}{\partial v_2} & = & -y,\\
\frac{\partial c_2}{\partial v_1} & = & -\tau,& \quad & \frac{\partial c_2}{\partial v_2} & = & 0.
\end{array}
$$
By Lemma \ref{2ff} we can compute:
\begin{eqnarray*}
\bar{\nabla}_{V_1}V_1 & = &\left(x_{11}+\frac{\lambda}{2}\nu y_1\right)X+\left(y_{11}-\frac{\lambda}{2}\nu x_1\right)Y-\frac{k\tau}{2}Z,\\
\bar{\nabla}_{V_2}V_2 & = &\left(\frac{\lambda}{2}(2c-r^2)-1\right)\left(xX+yY\right),\\
\bar{\nabla}_{V_1}V_2 & = & \bar{\nabla}_{V_2}V_1\ =\ \left(-y_1+\frac{\lambda}{4}\left(\nu x+(2c-r^2)y_1\right)\right)X\\
&&\phantom{aaaaaaaaaa}+\left(x_1+\frac{\lambda}{4}\left(\nu y-(2c-r^2)x_1\right)\right)Y-\frac{\tau}{2}Z.
\end{eqnarray*}
Since the unit normal is 
$$
\nuv= \frac{\lambda\left(y_1(2c-r^2)-x\nu\right)X-\lambda\left(x_1(2c-r^2)+y\nu\right)Y+2\tau Z}{\sql\sqrt{4\tau^2+{\lambda}\left(4c^2-4c\tau^2+\tau^4+\tau^2\nu^2\right)}},
$$
the components of the second fundamental form of $\mm$ with respect to $V_1$ and $V_2$ are:
\begin{eqnarray*}
h_{11} & = & \frac{\sql\left(-2ck+\frac{\lambda}{2}\nu\left(2c-\tau^2\right)\right)}{\sqrt{4\tau^2+\lambda\left(4c^2-4c\tau^2+\tau^4+\tau^2\nu^2\right)}}\\
h_{22} & = & -\frac{c\sql\nu\left(\lambda\left(2c-r^2\right)-2\right)}{{\sqrt{4\tau^2+\lambda\left(4c^2-4c\tau^2+\tau^4+\tau^2\nu^2\right)}}}\\
h_{12} & = & h_{21}\ = \ \frac{\sql\left({\lambda}\left((2c-r^2)^2-\tau^4-\tau^2\nu^2\right)-8c\right)}{4{\sqrt{4\tau^2+\lambda\left(4c^2-4c\tau^2+\tau^4+\tau^2\nu^2\right)}}}.
\end{eqnarray*}
Computing $h_{11}$ we used the fact that $\gamma$ is parametrized by arc length, hence $\left\langle \gamma',\gamma''\right\rangle=0$ and 
$$
\left\langle \gamma,\gamma''\right\rangle=\tau\left\langle T,\gamma''\right\rangle+\nu\left\langle N,\gamma''\right\rangle=k\nu.
$$
Since the inverse of the induced metric is
$$
g^{-1}=\frac{4}{4\tau^2+\lambda\left(4c^2-4c\tau^2+\tau^4+\tau^2\nu^2\right)}\left(\begin{array}{cc}
 r^2+\frac{\lambda}{4}\left(2c-r^2\right)^2& -\frac{\nu}{4}\left(\lambda\left(2c-r^2\right)-4\right)\\
-\frac{\nu}{4}\left(\lambda\left(2c-r^2\right)-4\right) & 1+\frac{\lambda}{4}\nu^2
\end{array}\right).
$$
The formula for $H=tr(Ag^{-1})$ can be derived after some standard computations. Finally, imposing the vertical translating equation \eqref{vert trasl}
$$
H=\bar g_{\lambda}\left(\nuv,\frac{Z}{\sql}\right)=\frac{2\tau}{\sqrt{4\tau^2+{\lambda}\left(4c^2-4c\tau^2+\tau^4+\tau^2\nu^2\right)}}
$$
we can derive the formula for the curvature $k$ after some algebraic manipulations. 
\cvd

The qualitative behavior of the curve $\gamma$ generating an helicoidal vertical translator is, for every $c$ and every $\lambda$, analogous to the solution of the same problem in the Euclidean space (Theorem 4.1 of \cite{Ha})

\begin{Theorem}\label{teor rototrasl}
Then there is a one-parameter family of helicoidal vertical translator: the generating curve $\gamma$ solution of \eqref{tau-nu} with $k$ as in \eqref{eq k} has exactly one point closest to the origin and consist of two properly embedded arms coming from this point which strictly go away to infinity and spiral infinitely many circles around the origin.
\end{Theorem}

Despite the general strategy of the proof is the same of \cite{Ha}, some new estimates are needed since in our case the expression of the curvature $k$ \eqref{eq k} is more complicated than the one found for the Euclidean space. One of the main difficulties appears when we try to count the numebers of zeros of $k$: with reference to the proof of Lemma 4.1 of \cite{Ha}, in our case, if $k=0$, $k'$ does not necessary have a sign. We postpone this problem to Lemma \ref{l03}. On the other hand we simplify some parts of the strategy of \cite{Ha} proving more directly that $r$ is unbounded. In order to start the analysis of the curve satisfying equation \eqref{eq k}, we can note that it exists and is complete, i.e. defined on all $\mathbb R$,  as an application of Lemma 3.1 of \cite{Ha}. 

\begin{Lemma}\label{l01}
There are no equilibrium points.
\end{Lemma}
\proof
Suppose that an equilibrium point exists, then $\tau'=\nu'=0$ at this point. From \eqref{tau-nu} it follows that $k=0$ or $\tau=0$. Using \eqref{tau-nu} again, we can exclude $k=0$ because otherwise we have
$$
0=\tau'=1+k\nu=1.
$$
Then we can consider $\tau=0$. In this case we have another contradiction because
$$
0=1+k\nu=\frac{4\lambda}{4\nu^2+\lambda(2c-\nu^2)^2}>0.
$$
\cvd
\begin{Lemma}\label{l01.2}
The functions $\tau$ and $\nu$ have at most a finite number of zeros.
\end{Lemma}
\proof
By \eqref{eq k} we have that, whenever $\tau=0$
$$
\tau'=\frac{4\lambda c^2}{4\nu^2+\lambda(2c-\nu^2)^2}>0
$$
holds. This implies that $\tau$ has at most one zero. In a similar way, if $\nu=0$, from \eqref{eq k} we get that
$$
\nu'=\frac{\tau^2}{\sql c}\geq0.
$$
Then $\nu$ has at most one zero if $(0,0)\notin\gamma$ (i.e. $\tau\neq 0$) and at most three zeros if the curve passes from the origin.
\cvd

\begin{Lemma}\label{l02}
The function $r^2$ has exactly a global minimum and $\lim_{s\rightarrow\pm\infty}r^2=\infty$. In particular $\gamma$ has two embedded infinite arms starting from the minimum of $r^2$.
\end{Lemma}
\proof
Since $r^2=\nu^2+\tau^2$, by \eqref{tau-nu} we have that $(r^2)'=2\tau$. In view of Lemma \ref{l01.2}, we claim that $\tau$ has exactly a zero and therefore $r^2$ has a global minimum. In fact, if it is not the case, up to reverse the orientation of $\gamma$, we can suppose that $\tau>0$. Then $r^2$ is monotone and, in particular, there exists a constant $a\geq 0$ such that
$$
\lim_{s\rightarrow-\infty}r^2=a.
$$
Therefore $\lim_{s\rightarrow-\infty}(r^2)'=0$, hence
\begin{equation*}
\lim_{s\rightarrow-\infty}\tau=0,\quad
\lim_{s\rightarrow-\infty}\nu^2=a.
\end{equation*}
Using these equality and \eqref{eq k}, similarly to Lemma \ref{l01}, we get
$$
0=\lim_{s\rightarrow-\infty}\tau'=\lim_{s\rightarrow-\infty}1+k\nu=\frac{4\lambda c^2}{4c+\lambda(2c-a)^2}>0,
$$
giving a contradiction.
We proved that $r^2$ has exactly one critical point, then it has a limit for $s\rightarrow\pm\infty$. Repeating the argument just seen, we can exclude the convergence to a finite limit.
\cvd

\begin{Lemma}\label{l03}
There exist a compact interval $I\subset \mathbb R$ containing all the zeros of $k$.
\end{Lemma}
\proof
Fix a point where $k=0$, by \eqref{tau-nu} we have $\tau'=1$ and $\nu'=0$ in this point. Deriving \eqref{eq k} in this point we have:
\begin{equation}\label{eq k'}
k'=\frac{-3\lambda\tau^2\nu^2-4\lambda^{\frac 32}c\tau\nu-5\lambda\tau^4+12(\lambda c-1)\tau^2-4\lambda c^2}{\sql c\left(4r^2+\lambda(2c-r^2)^2\right)}.
\end{equation}
Note that, contrary to the Euclidean case \cite{Ha}, in general $k'$ does not have a sign. We can study the set where $k$ and $k'$ vanishes as curves in the real projective plane with homogeneous coordinates $(\tau,\nu,\zeta)$ depending on the two paramethers $\lambda$ and $c$. Let us define the following polynomials:
\begin{eqnarray*}
P_0(\tau,\nu,\zeta) & = & \lambda^{\frac 32}c\nu\zeta^2(4c\zeta^2-2\tau^2-\nu^2)-4\sql c\nu\zeta^4\\
&&-\tau(4\tau^2\zeta^2+\lambda(4c^2\zeta^4-4c\tau^2\zeta^2+\tau^4+\tau^2\nu^2));\\
P_1(\tau,\nu,\zeta) & = & -3\lambda\tau^2\nu^2-4\lambda^{\frac 32}c\tau\nu\zeta^2-5\lambda\tau^4+12(\lambda c-1)\tau^2\zeta^2-4\lambda c^2\zeta^4,
\end{eqnarray*}
and let $\mathcal C_0=\left\{P_0=0\right\}$ be the curve where $k$ vanishes and $\mathcal C_1=\left\{P_1=0\right\}$ the curve where $k'=0$. The point $q=(0,1,0)$ is the only point at infinity of $\mathcal C_0$. Note that $q$ belongs to $\mathcal C_1$ too. Deriving $P_0$ three times we can see that $\mathcal C_0$ is asymptotic to the curve $\mathcal A$ of equation $\sql c \nu\zeta^2+\tau^3=0$. Since in Lemma \ref{l02} we proved that $r^2=\tau^2+\nu^2$ diverges,  if $|s|$ is big enough, in order to simplify computations, we can consider points on $\mathcal A$ rather than on $\mathcal C_0$. We have that
$$
P_1\left(\tau,-\frac{\tau^3}{\sql c},1\right)=-\frac{3}{c^2}\tau^8-\lambda\tau^4+12(\lambda c-1)\tau^2-4\lambda c^2<-1
$$
holds if $\tau^2$ is big enough. Increasing $|s|$ if necessary, we have that if $k(s)=0$ then $k'(s)<0$, therefore outside a compact set the curvature $k$ cannot be zero.\cvd

\begin{Lemma}\label{l04}
Both $\tau$ and $\nu$ have limits when $s\rightarrow\pm\infty$. In particular $\tau$ has finite limit in each direction, while $\lim_{s\rightarrow\pm\infty}\nu=\pm\infty$.
\end{Lemma}
\proof
The proof for the function $\tau$ uses a strategy similar to that of Lemma \ref{l03}. Using the explicit expression of $k$ \eqref{eq k}, passing to homogeneous coordinates like in the proof of Lemma \ref{l03}, we have that $\tau'=1+k\nu=0$ if and only if we are on the curve of equation
$$
P(\tau,\nu,\zeta)=\lambda^{\frac 32}c\left(2c\zeta^2-\tau^2\right)^2\zeta^2-\lambda\tau\nu\left(\left(2c\zeta^2-\tau^2\right)^2+\tau^2\nu^2\right)+4\sql c\tau^2\zeta^4-4\tau^3\nu\zeta^2=0.
$$
This curve has two points at infinity: $q_1=(1,0,0)$ and $q_2=(0,1,0)$ and it is asymptotic to the lines given by $\tau\nu=0$. This means that if $\tau'(s)=0$ for an $|s|$ large enough, then one of $\tau$ or $\nu$ is close to zero (but not both because of Lemma \ref{l02}). Now we want to estimate $\tau''$ in these points. Let us define by $F$ (resp. $G$) the numerator (resp. the denominator) of $k$ in \eqref{eq k}, using \eqref{tau-nu} and the fact that in critical points of $\tau$ $k=-\frac{1}{\nu}$, we get:
\begin{eqnarray*}
\tau''&=&k'\nu-k^2\tau\ =\  \frac{F'\nu^3+G'\nu^2-\tau G}{\nu^2 G}.
\end{eqnarray*}
Since $G$ is positive, the sign of $\tau''$ is given by $F'\nu^3+G'\nu^2-\tau G$. Recalling that in the points considered we have $\tau'=0$ and $\nu'=\frac{\tau}{\nu}$, after some computations we have that:
\begin{equation}\label{FG}
F'\nu^3+G'\nu^2-\tau G=\tau\left(2\lambda\tau^3\nu^2\left(\sql c-\nu\tau\right)-4\sql c\tau^2-\lambda^{\frac 32}c\left(2c-\tau^2\right)^2\right).
\end{equation}
Suppose that we have $|s|$ big enough and $\nu$ close to zero, in this case from \eqref{FG} we have that
$$
\tau''\nu^2G\approx -\sql c\tau\left(4\tau^2+\lambda(2c-\tau^2)^2\right),
$$
hence the sign of $\tau''$ is determined by the sign of $\tau$: $\tau$ has a local maximum (resp. minimum) if $\tau>0$ (resp. $\tau<0$). This means that, outside a compact set of $\mathbb R$, $\tau$ does not have critical points with $\nu$ too small.

Suppose now that $|s|$ is big and $\tau$ is close to zero. By \eqref{FG} we get
$$
\tau''\nu^2 G\approx-4\tau\lambda^{\frac 32}c^3.
$$
Therefore once again we can deduce the sign $\tau''$ from that of $\tau$. Arguing as above we can find a compact set that encloses all the possible critical points of $\tau$. Outside this compact set $\tau$ is monotone, hence we have that $\tau$ has a limit in each directions. Since in Lemma \ref{l01.2} we proved that $\tau$ has exactly one zero, we get that the two limits are finite. 

For the function $\nu$ the proof is simpler. Since $\nu'=-k\tau$, enlarging the interval $I$ in order to include the zero of $\tau$ too, by Lemma \ref{l01.2} and Lemma \ref{l03}, we have that $\nu$ does not have critical points outside $I$. Therefore this function too has a limit in each direction. Moreover since $\nu^2=r^2-\tau^2$, by Lemma \ref{l02} and what we showed so far in this proof, we have that $\nu^2$ diverges in each direction. In order to determinate the sign of $\nu$, we recall that in Lemma \ref{l01.2} we proved that $\nu$ has a finite numbers of zeros and $\nu'>0$ in these points.
\cvd

Finally we can conclude the proof of the Theorem \ref{teor rototrasl} with the following series of results. We omit the proof because it is the same of Lemmas 4.5, 4.7, 4.8 and 4.9 of \cite{Ha}. We point out that we found opposite signs because we have the opposite sign in the translator equation \eqref{vert trasl}.

\begin{Lemma}
\begin{itemize}
\item[1)] $\lim_{s\rightarrow\pm\infty}k=0^{\mp}$;
\item[2)] each arm of $\gamma$ has infinite total curvature, i.e. $\lim_{s\rightarrow\pm\infty}\vartheta=-\infty$;
\item[3)] the limit growing direction of the arms is given by $\lim_{s\rightarrow\pm\infty}\frac{rT}{\gamma}=\mp i$;
\item[4)] writing $\gamma=re^{i\vartheta}$ for some angle $\vartheta$, we have that $\lim_{s\rightarrow\pm\infty}\vartheta=-\infty$, i.e. each arm spirals infinitely many times around the origin.
\end{itemize}
\end{Lemma}

We finish this section showing that there are no other helicoidal vertical translators.

\begin{Theorem}
For every $(a,b,c)\in\mathbb R^3$ with $(a,b)\neq(0,0)$ there are no vertical translators invariant under the action of the group $G=\left\{L_{(au,bu,cu)}\circ\rho_u\left|u\in\mathbb R\right.\right\}$.
\end{Theorem}
\proof Let $\mm$ be a $G$-invariant surface. As seen before, to horizontal rotations we can suppose that $b=0$, therefore there is a planar curve $\gamma=(\gamma_1,\gamma_2)$ such that 
\begin{eqnarray*}
\mm(v_1,v_2) & = & \left(\gamma_1(v_1)\cos(v_2)-\gamma_2(v_1)\sin(v_2)+av_2,\phantom{\frac{av_2}{2}} \right.\\
&&\phantom{a.}\gamma_1(v_1)\sin(v_2)+\gamma_2(v_1)\cos(v_2), \\
&&\phantom{a.}\left.cv_2+\frac{av_2}{2}\left(\gamma_1(v_1)\sin(v_2)+\gamma_2(v_1)\cos(v_2)\right)\right)
\end{eqnarray*}
By Lemma \ref{2ff} a tangent basis of $\mm$ is 
$$
V_1=x_1X+y_1Y+c_1Z,\quad V_2=(a-y)X+(x-av_2)Y+c_2Z;
$$
where \begin{equation*}
c_1  =  \frac{1}{2}\left(\nu+ay_1v_2\right),\quad c_2  =  \frac 12\left(2b-r^2+2ay+2axv_2-a^2v_2^2\right).
\end{equation*}
Let us define the following quantities:
\begin{equation*}
\begin{array}{rcl}
E_1 & = & \det\left(\begin{array}{cc}
c_1 & y_1\\
c_2 & x-av_2
\end{array}\right)= -\frac{a}{2}\left(\nu+xy_1\right)v_2+\frac 12\left(\nu x-2y_1b^2+2y_1r^2-2ayy_1\right);\\
\\
E_2 & = & \det\left(\begin{array}{cc}
x_1 & c_1\\
a-y & c_2
\end{array}\right)\\
& = & -\frac{a^2x_1}{2}v_2^2+\frac{a}{2}\left(2xx_1-ay_1+yy_1\right)v_2+\frac 12\left(2b^2x_1-r^2x_1+2ax_1y-a\nu+\nu y\right);\\
\\
E_3 & = & \tau-ay_1-ax_1v_2.
\end{array}
\end{equation*}
A unit normal vector field is 
\begin{equation}\label{norm-heli}
\nu= \frac{-\lambda E_1X-\lambda E_2Y+E_3Z}{\sqrt{\lambda^2 E_1^2+\lambda^2 E_2^2+\lambda E_3^2}}.
\end{equation}
If $\mm$ is a vertical translator, because of the $G$-invariance, the function
$$
\frac{\tau-a(y_1+v_2x_1)}{\sqrt{\lambda E_1^2+\lambda E_2^2+E_3^2}}=\bar g_{\lambda}\left(\nuv,\frac{Z}{\sql}\right)=H
$$
has to be independent on $v_2$, but this is impossible if $a\neq 0$ because,in this case, this property holds if only if $\gamma_1'=\gamma_2'=0$. The proof in a consequence of the fact that the functions $1,\ \cos(v_2),\ \sin(v_2),\ v_2\cos(v_2),\ v_2\sin(v_2),\ v_2^2\sin(v_2),\ v_2^2\cos(v_2)$ are linearly independent.
\cvd

\section{Translator with respect to a generic direction}
In this final section we wish to study the invariant surfaces translating in a fixed non-vertical direction proving Theorem \ref{main 3}. Fix a the Killing vector field $$V=a_1F_1+a_2F_2+a_3F_3=a_1X+a_2Y+(a_3+a_1y-a_2x)Z,$$ with $(a_1,a_2,a_3)\in\mathbb R^3$ and $(a_1,a_2)\neq(0,0)$. First of all we note that in this case it is not possible to require a translation with unit speed because the norm of $V$ is not constant. 
\\

\noindent\emph{Proof of Theorem \ref{main 3}}. 
\begin{itemize}
\item[1)] If $\mm$ is invariant for vertical translation, we can find a planar curve $\gamma$ such that $\mm=\gamma\times\mathbb R$. Since the natural projection $\pi:(x,y,z)\in\hhh_3\mapsto(x,y)\in\mathbb R^2$ is a Riemannian submersion with geodesic fibers, we have that the mean curvature of $\mm$ in a point $p$ is equal to the curvature of $\gamma$ in $\pi(p)$ and the mean curvature flow commutes with $\pi$, as proven in \cite{Pi}. Moreover $Z$ is tangent to $\mm$ hence the translator equation becomes
$$
k=H=\bar g_{\lambda}(\nuv,V)=\left\langle \nuv,(a_1,a_2)\right\rangle,
$$
where, with an abuse of notation, we used the same symbol for the normal of $\mm$ in $\hhh_3$ and the normal of $\gamma$ in $\mathbb R^2$. Hence $\gamma$ satisfies the equation for a translator in the Euclidean plane and it is well known that the grim reaper is the only solution.
\item[2)] Let $\mm$ be invariant by the action of the group generated by $L_{(x_0,y_0,z_0)}$. As seen many times in this paper, it is not restrictive to choose $x_0=1$ and $y_0=0$. In this case the Killing vector field $F_1+z_0F_3$ is tangent to $\mm$. Since trivially $V=a_1(F_1+z_0F_3)+a_2F_2+(a_3-a_1z_0)F_3$, then
\begin{eqnarray*}
\bar g_{\lambda}(\nuv,V)=\bar g_{\lambda}(\nuv,a_2F_2+(a_3-a_1z_0)F_3).
\end{eqnarray*}
Moreover from \eqref{mc graph lambda} we have that the mean curvature of $H$ does not depend on $x$, then $a_2=0$. Then the thesis follows choosing $c=a_1$.
\item[3)] It is easy to prove that the round cylinders with vertical axis are not translators: they shrink to the $z$-axis in finite time. Therefore if $\mm$ is invariant by horizontal rotation we have that at least locally, it can be written as an horizontal graph $z=\varphi(r)$. Its mean curvature is given by \eqref{H rot}, in particular it depends only on $r$. On the other hand, by Lemma \ref{2ff} we can compute the normal vector, hence 
\begin{eqnarray*}
\bar g_{\lambda}(\nuv,V) & = & \frac{\sql}{\sqrt{1+\lambda\left((\varphi')^2+\frac{r^2}{4}\right)}}\left(-x\left(a_1\frac{\varphi'}{r}+\frac{a_2}{2}\right)-y\left(a_2\frac{\varphi'}{r}-\frac{a_1}{2}\right)+a_3\right).
\end{eqnarray*}
If we hope to find a solution for the translating equation, we need that also this quantity depends only on $r$ (and not on $x$ and $y$ independently). Therefore we require that
$$
a_1\frac{\varphi'}{r}+\frac{a_2}{2}=a_2\frac{\varphi'}{r}-\frac{a_1}{2}=0.
$$
As $(a_1,a_2)\neq(0,0)$, it is easy to deduce that this condition can be verified for every $r$ if and only if $\varphi'$ vanishes everywhere, but in this case the associated $\mm$ is a horizontal plane that is not a translator: this surface is minimal, hence it does not move, but $V$ does not belong to its tangent plane.
\item[4)] Let $\mm$ be a helicoid of general type, as in the previous section, up to horizontal rotation, it can be written as 
$$\mm(v_1,v_2)=L_{(av_2,0,cv_2)}\circ\rho_{v_2}(\gamma_1(v_1),\gamma_2(v_1),0).$$ 
The unit normal vector of $\mm$ is given by \eqref{norm-heli}. Once again the function $\bar g_{\lambda}(\nuv,V)$ has to share the same symmetries of $\mm$, in particular it does not depend on $v_2$. Suppose that there is a function $F=F(v_1)$ such that $\bar g_{\lambda}(\nuv,V)=F(v_1)$. Using \eqref{norm-heli} and the notation introduced in the previous section we get that it is equivalent to
\begin{equation}\label{EEE}
\left(a_1\sql E_1+a_2\sql E_2-\frac{a_3}{\sql}E_3\right)^2=F^2\left(\lambda E_1^2+\lambda E_2^2+E_3^2\right) 
\end{equation}
for every $v_1,\ v_2$. Roughly speaking this property holds if the coefficients of all the terms in $v_2$ vanish. For example we can test this equation on the sequence $v_2=\frac{\pi}{2}m$, with $m\in\mathbb Z$ and we want that \eqref{EEE} does not depend on $m$. Since $a_1^2+a_2^2\neq 0$ and $c\neq 0$ we derive that it is possible if and only if $\gamma_1'=\gamma_2'=0$, having a contradiction.

\cvd
\end{itemize}

\begin{Remark}
Theorem \ref{main 3}, part $4)$ says in particular that in $\hhh_3$ the only helicoidal surfaces that are translators are those with vertical axis, i.e. when $a=b=0$. Obviously in the Euclidean space we can find helicoidal translator with any axis. In fact such an helicoid is always isometric to a helicoid with vertical axis: it is sufficient to rotate one axis in the other. This is no more true in our setting because the lack of symmetries in $\hhh_3$.
\end{Remark}

\bigskip

\noindent Giuseppe Pipoli, \emph{Department of Information Engineering, Computer Science and Mathematics, Universit\`a degli Studi dell'Aquila}, via Vetoio 1, 67100, L'Aquila, Italy.\\
Email: giuseppe.pipoli@univaq.it
\end{document}